\documentclass{amsart} 

\usepackage{amssymb} 
\usepackage{amsmath} 
\usepackage{amscd} 
\usepackage{latexsym} 
\usepackage{amsthm} 
\usepackage{verbatim} 
\usepackage{graphics} 
\usepackage{epsfig}

\let\cal\mathcal 
\def\Ascr{{\cal A}} 
\def\Bscr{{\cal B}} 
\def\Cscr{{\cal C}} 
\def\Dscr{{\cal D}} 
\def\Escr{{\cal E}} 
\def\Fscr{{\cal F}} 
\def\Gscr{{\cal G}} 
\def\Hscr{{\cal H}} 
\def\Iscr{{\cal I}} 
\def\Jscr{{\cal J}} 
\def\Kscr{{\cal K}} 
\def\Lscr{{\cal L}} 
\def\Mscr{{\cal M}} 
\def\Nscr{{\cal N}} 
\def\Oscr{{\cal O}} 
\def\Pscr{{\cal P}} 
\def\Qscr{{\cal Q}} 
\def\Rscr{{\cal R}} 
\def\Sscr{{\cal S}} 
\def\Tscr{{\cal T}} 
\def\Uscr{{\cal U}} 
\def\Vscr{{\cal V}} 
\def\Wscr{{\cal W}} 
\def\Xscr{{\cal X}} 
\def\Yscr{{\cal Y}} 
\def\Zscr{{\cal Z}} 

\let\blb\mathbb 
\def\AA{{\blb A}} 
\def\BB{{\blb B}} 
\def\CC{{\blb C}}
\def\DD{{\blb D}}  
\def\FF{{\blb F}} 
\def\GG{{\blb G}} 
\def \HH{{\blb H}} 
\def\II{{\blb I}}
\def\LL{{\blb L}} 
\def \NN{{\blb N}} 
\def \PP{{\blb P}} 
\def\QQ{{\blb Q}} 
\def \RR{{\blb R}}
\def\VV{{\blb V}}
\def\WW{{\blb W}}
\def\UU{{\blb U}}
\def \TT{{\blb T}}
\def\XX{{\blb X}}
\def \ZZ{{\blb Z}}

\def\1{{\mathbb{1}}}

\def\a{{\tilde{\alpha}}} 
\def\b{{\tilde{\beta}}} 
\def\c{{\tilde{\gamma}}} 

\let\st\ast 
\let\at\ast 
\def\Ab{\mathbb{Ab}} 
\def\add{\operatorname {add}} 
\def\adm{\operatorname{adm}} 
\def\Ann{\operatorname{Ann}} 
\def\Bimod{\operatorname{Bimod}} 
\def\cd{\operatorname {cd}} 
\def\ch{\mathop{\text{Ch}}} 
\def\Ch{\mathop{\mathrm{Ch}}} 
\def\CH{\mathop{\mathrm{CH}}} 
\def\charact{\operatorname{char}} 
\def\codim{\operatorname {codim}} 
\def\coh{\mathop{\text{\upshape{coh}}}} 
\def\coker{\operatorname {coker}} 
\def\cone{\operatorname {cone}} 
\def\coop{{\operatorname {coop}}} 
\def\ctimes{\mathbin{\hat{\otimes}}} 
\def\ctensor{\mathbin{\hat{\otimes}}} 
\def\d{\downarrow} 
\def\D{\operatorname{D}} 
\def\depth{\operatorname {depth}} 
\def\Der{\operatorname{Der}} 
\def\diag{\operatorname {diag}} 
\def\dis{\operatorname {dis}} 
\def\Dis{\operatorname{Dis}} 
\def\End{\operatorname {End}} 
\def\exist{\exists} 
\def\Ext{\operatorname {Ext}} 
\def\fid{\operatorname{fid}} 
\def\Filt{\operatorname {Filt}} 
\def\for{\operatorname {For}} 
\def\For{\operatorname {For}} 
\def\fpd{\operatorname{fpd}} 
\def\Fract{\operatorname {Fract}} 
\def\G{\mathop{\underline{\underline{{\Gamma}}}}\nolimits} 
\def\Gal{\operatorname {Gal}} 
\def\gk{\operatorname {gk}} 
\def\gkdim{\operatorname {gkdim}} 
\def\Gl{\operatorname {Gl}} 
\def\GL{\operatorname {GL}} 
\def\gldim{\operatorname {gl\,dim}} 
\def\gr{\operatorname{gr}} 
\def\Gr{\operatorname{Gr}} 
\def\GR{\operatorname{GR}} 
\def\grmod{\operatorname {grmod}} 
\def\GrMod{\operatorname {GrMod}} 
\def\Hilb{\operatorname {Hilb}} 
\def\Hol{\text{-Hol}} 
\def\Hom{\operatorname {Hom}} 
\def\ICh{\mathop{\text{ICh}}} 
\def\id{{\operatorname {id}}} 
\def\Id{\operatorname{id}} 
\def\im{\operatorname {im}} 
\def\inv{\operatorname {inv}}
\def\K0{{K_{0}(\blb P^{2}_{q})}}
\def\kard{\operatorname {kard}}  
\def\ker{\operatorname {ker}} 
\def\Ker{\operatorname {ker}} 
\def\l{\leftarrow} 
\def\la{\operatorname{\lambda}} 
\def\length{\mathop{\text{length}}} 
\def\Lie{\mathop{\text{Lie}}} 
\def\Lotimes{\overset{{\bf L}}{\otimes}}
\def\LotimesB{\overset{{\bf L}}{\otimes}_{B}\hspace{1mm}}  
\def\LotimesD{\overset{{\bf L}}{\otimes}_{D}}
\def\LotimesDelta{\overset{{\bf L}}{\otimes}_{\Delta}\hspace{1mm}}
\def\lr{\longrightarrow} 
\def\mod{\operatorname{mod}} 
\def\Mod{\operatorname{Mod}} 
\def\Nil{\operatorname{Nil}} 
\def\opp{{\operatorname{opp}}} 
\def\ord{\operatorname {ord}} 
\def\P2q{\operatorname {{\blb P}^{2}_{q}}}  
\def\pc{\operatorname {pc}} 
\def\PC{\operatorname {PC}} 
\def\pd{\operatorname {pd}} 
\def\pdim{\operatorname {pd}} 
\def\PGL{\operatorname {PGL}} 
\def\Pic{\operatorname {Pic}} 
\def\PP{\operatorname {\blb P}} 
\def\pr{\mathop{\text{pr}}\nolimits} 
\def\pwqch{\text{-(w)qch}} 
\def\qch{\text{-qch}} 
\def\Qch{\operatorname{Qch}} 
\def\Qcoh{\operatorname {Qcoh}} 
\def\qgr{\operatorname{qgr}} 
\def\QGr{\operatorname{QGr}} 
\def\qpc{\operatorname {qpc}} 
\def\QSch{\mathop{\text{Qsch}}} 
\def\quot{/\!\!/} 
\def\r{\rightarrow} 
\def\rad{\operatorname {rad}} 
\def\rank{\operatorname {rank}} 
\def\red{\operatorname {red}}
\def\Rees{\operatorname {Rees}}
\def\Rep{\operatorname {Rep}} 
\def\relint{\operatorname {relint}} 
\def\Res{\operatorname{Res}} 
\def\RHom{\operatorname {RHom}} 
\def\rk{\operatorname {rk}} 
\def\sh{\operatorname{sh}} 
\def\sign{\operatorname{sign}} 
\def\Sl{\operatorname {Sl}} 
\def\Skl{\operatorname{Skl}} 
\def\Span{\operatorname {Span}}
\def\spec{\operatorname {Spec}} 
\def\Spec{\operatorname {Spec}} 
\def\supp{\operatorname {supp}} 
\def\Supp{\mathop{\text{\upshape Supp}}} 
\def\Tails{\operatorname {Tails}}  
\def\tails{\operatorname {tails}}  
\def\Tau{\mathcal{T}} 
\def\Tor{\operatorname {Tor}} 
\def\tot{\operatorname {Tot}} 
\def\Tot{\operatorname {Tot}} 
\def\tr{\operatorname {Tr}} 
\def\Tr{\operatorname {Tr}} 
\def\u{\uparrow} 
\def\un{\operatorname {--}}  
\def\uRHom{\operatorname {R\mathcal{H}\mathit{om}}} 
\def\wqch{\text{-wqch}} 
\def\weight{\operatorname{weight}} 

\DeclareMathOperator{\Proj}{Proj} 
\DeclareMathOperator{\Pro}{Pro} 
\DeclareMathOperator{\Ind}{Ind} 
\DeclareMathOperator{\Alg}{Alg} 
\DeclareMathOperator{\GrAlg}{GrAlg} 
\DeclareMathOperator{\Rel}{Rel} 
\DeclareMathOperator{\Tors}{Tors} 
\DeclareMathOperator{\tors}{tors} 
\DeclareMathOperator{\HTor}{\mathcal{T}\mathit{or}} 
\DeclareMathOperator{\HHom}{\mathcal{H}\mathit{om}} 
\DeclareMathOperator{\HEnd}{\mathcal{E}\mathit{nd}} 
\DeclareMathOperator{\Nrd}{Nrd} 
\DeclareMathOperator{\Aut}{Aut} 
\DeclareMathOperator{\ram}{ram} 
\DeclareMathOperator{\Br}{Br} 
\DeclareMathOperator{\Div}{Div} 
\DeclareMathOperator{\Prin}{Prin} 
\DeclareMathOperator{\cor}{cor} 
\DeclareMathOperator{\res}{res} 

\let\dirlim\injlim 
\let\invlim\projlim 
\newtheorem*{theoremA}{Theorem A} 
\newtheorem*{theoremB}{Theorem B}

\newtheorem{lemma}{Lemma}[section] 
\newtheorem{proposition}[lemma]{Proposition}

\theoremstyle{definition} 

\newtheorem{example}[lemma]{Example}

{

\newtheorem{case}{Case}

\newtheorem{problem}{Problem}
\newtheorem{solution}{Solution}
} 

\theoremstyle{remark} 

\newtheorem{remark}[lemma]{Remark} 
 
\def\thenotation{} 

\newdimen\uboxsep \uboxsep=1ex 
\def\uboxn#1{\vtop to 0pt{\hrule height 0pt depth 0pt\vskip\uboxsep 
\hbox to 0pt{\hss #1\hss}\vss}} 

\def\uboxs#1{\vbox to 0pt{\vss\hbox to 0pt{\hss #1\hss} 
\vskip\uboxsep\hrule height 0pt depth 0pt}}

\def\debug{\tracingall} 
\def\enddebug{\showlists} 

\numberwithin{equation}{section} 

\def\a{{\tilde{\alpha}}} 
\def\b{{\tilde{\beta}}} 
\def\c{{\tilde{\gamma}}} 

\def\debug{\tracingall} 
\def\enddebug{\showlists} 

\keywords{Partitions of integers, Ferrers graphs} 
\subjclass{Primary 05A17, 11P81} 
\author{Koen De Naeghel and Nicolas Marconnet} 

\address{Koen De Naeghel \\Departement WNI\\Hasselt University\\ Agoralaan gebouw D \\ B-3590 
Diepenbeek (Belgium).}
\email[K. De Naeghel]{koen.denaeghel@uhasselt.be}
\address{Nicolas Marconnet \\ 
Departement Wiskunde en Informatica\\Universiteit Antwerpen\\Middelheimlaan 1\\ B-2020 Antwerp (Belgium).}
\email[N. Marconnet]{%nicolas.marconnet@univ-st-etienne.fr
nicolas.marconnet@ua.ac.be} 

\date{January 5, 2006} 
\title[An inequality on broken chessboards]{An inequality on broken chessboards} 

\begin{document} 

\begin{abstract} 
For any partition of a positive integer we consider the chess (or draughts) colouring of its associated Ferrers graph. Let $b$ denote the total number of black unit squares, and $w$ the number of white squares. In this note we characterize all pairs $(b,w)$ which arise in this way. This simple combinatorical result was discovered by characterizing Hilbert series of certain right modules over cubic three-dimensional Artin-Schelter algebras. However in this note we present a purely combinatorical proof. \\

The result is (at least partially) known in literature \cite[Problem 10]{competition}, however we found it interesting to present an elementary proof. All additional references and remarks will be mostly appreciated.
\end{abstract} 

\maketitle 

\tableofcontents 

\section{Introduction}

A partition of a positive integer $n$ is a finite nonincreasing sequence of positive integers $\la_{1},\la_{2}, \dots, \la_{r}$ such that $\sum_{i = 1}^{r}\la_{i} = n$. We denote $\la = (\la_{1},\la_{2}, \dots, \la_{r})$. To each partition $\la$ is associated its Ferrers graph: a pattern of unit squares with the $i$-th row (counting from $i = 0$) having $la_{i+1}$ unit squares (see \S\ref{2.1} for a more formal definition). As an example the Ferrers graph of the partition $\la = (8,6,6,5,2,1,1)$ of $29$ is given by \\

%\input partition29.tex \\
%Created by jPicEdt 1.x
%Standard LaTeX format (emulated lines)
%Wed Jul 06 22:51:42 CEST 2005
\unitlength 1mm
\begin{picture}(80.00,35.00)(0,0)

\linethickness{0.15mm}
%Rectangle(40.00,0.00)(45.00,5.00)  
\put(40.00,0.00){\line(1,0){5.00}}
\put(40.00,0.00){\line(0,1){5.00}}
\put(45.00,0.00){\line(0,1){5.00}}
\put(40.00,5.00){\line(1,0){5.00}}
%End Rectangle

\linethickness{0.15mm}
%Rectangle(45.00,0.00)(50.00,5.00)  
\put(45.00,0.00){\line(1,0){5.00}}
\put(45.00,0.00){\line(0,1){5.00}}
\put(50.00,0.00){\line(0,1){5.00}}
\put(45.00,5.00){\line(1,0){5.00}}
%End Rectangle

\linethickness{0.15mm}
%Rectangle(50.00,0.00)(55.00,5.00)  
\put(50.00,0.00){\line(1,0){5.00}}
\put(50.00,0.00){\line(0,1){5.00}}
\put(55.00,0.00){\line(0,1){5.00}}
\put(50.00,5.00){\line(1,0){5.00}}
%End Rectangle

\linethickness{0.15mm}
%Rectangle(55.00,0.00)(60.00,5.00)  
\put(55.00,0.00){\line(1,0){5.00}}
\put(55.00,0.00){\line(0,1){5.00}}
\put(60.00,0.00){\line(0,1){5.00}}
\put(55.00,5.00){\line(1,0){5.00}}
%End Rectangle

\linethickness{0.15mm}
%Rectangle(60.00,0.00)(65.00,5.00)  
\put(60.00,0.00){\line(1,0){5.00}}
\put(60.00,0.00){\line(0,1){5.00}}
\put(65.00,0.00){\line(0,1){5.00}}
\put(60.00,5.00){\line(1,0){5.00}}
%End Rectangle

\linethickness{0.15mm}
%Rectangle(65.00,0.00)(70.00,5.00)  
\put(65.00,0.00){\line(1,0){5.00}}
\put(65.00,0.00){\line(0,1){5.00}}
\put(70.00,0.00){\line(0,1){5.00}}
\put(65.00,5.00){\line(1,0){5.00}}
%End Rectangle

\linethickness{0.15mm}
%Rectangle(70.00,0.00)(75.00,5.00)  
\put(70.00,0.00){\line(1,0){5.00}}
\put(70.00,0.00){\line(0,1){5.00}}
\put(75.00,0.00){\line(0,1){5.00}}
\put(70.00,5.00){\line(1,0){5.00}}
%End Rectangle

\linethickness{0.15mm}
%Rectangle(75.00,0.00)(80.00,5.00)  
\put(75.00,0.00){\line(1,0){5.00}}
\put(75.00,0.00){\line(0,1){5.00}}
\put(80.00,0.00){\line(0,1){5.00}}
\put(75.00,5.00){\line(1,0){5.00}}
%End Rectangle

\linethickness{0.15mm}
%Rectangle(40.00,5.00)(45.00,10.00)  
\put(40.00,5.00){\line(1,0){5.00}}
\put(40.00,5.00){\line(0,1){5.00}}
\put(45.00,5.00){\line(0,1){5.00}}
\put(40.00,10.00){\line(1,0){5.00}}
%End Rectangle

\linethickness{0.15mm}
%Rectangle(45.00,5.00)(50.00,10.00)  
\put(45.00,5.00){\line(1,0){5.00}}
\put(45.00,5.00){\line(0,1){5.00}}
\put(50.00,5.00){\line(0,1){5.00}}
\put(45.00,10.00){\line(1,0){5.00}}
%End Rectangle

\linethickness{0.15mm}
%Rectangle(50.00,5.00)(55.00,10.00)  
\put(50.00,5.00){\line(1,0){5.00}}
\put(50.00,5.00){\line(0,1){5.00}}
\put(55.00,5.00){\line(0,1){5.00}}
\put(50.00,10.00){\line(1,0){5.00}}
%End Rectangle

\linethickness{0.15mm}
%Rectangle(55.00,5.00)(60.00,10.00)  
\put(55.00,5.00){\line(1,0){5.00}}
\put(55.00,5.00){\line(0,1){5.00}}
\put(60.00,5.00){\line(0,1){5.00}}
\put(55.00,10.00){\line(1,0){5.00}}
%End Rectangle

\linethickness{0.15mm}
%Rectangle(60.00,5.00)(65.00,10.00)  
\put(60.00,5.00){\line(1,0){5.00}}
\put(60.00,5.00){\line(0,1){5.00}}
\put(65.00,5.00){\line(0,1){5.00}}
\put(60.00,10.00){\line(1,0){5.00}}
%End Rectangle

\linethickness{0.15mm}
%Rectangle(65.00,5.00)(70.00,10.00)  
\put(65.00,5.00){\line(1,0){5.00}}
\put(65.00,5.00){\line(0,1){5.00}}
\put(70.00,5.00){\line(0,1){5.00}}
\put(65.00,10.00){\line(1,0){5.00}}
%End Rectangle

\linethickness{0.15mm}
%Rectangle(40.00,10.00)(45.00,15.00)  
\put(40.00,10.00){\line(1,0){5.00}}
\put(40.00,10.00){\line(0,1){5.00}}
\put(45.00,10.00){\line(0,1){5.00}}
\put(40.00,15.00){\line(1,0){5.00}}
%End Rectangle

\linethickness{0.15mm}
%Rectangle(45.00,10.00)(50.00,15.00)  
\put(45.00,10.00){\line(1,0){5.00}}
\put(45.00,10.00){\line(0,1){5.00}}
\put(50.00,10.00){\line(0,1){5.00}}
\put(45.00,15.00){\line(1,0){5.00}}
%End Rectangle

\linethickness{0.15mm}
%Rectangle(50.00,10.00)(55.00,15.00)  
\put(50.00,10.00){\line(1,0){5.00}}
\put(50.00,10.00){\line(0,1){5.00}}
\put(55.00,10.00){\line(0,1){5.00}}
\put(50.00,15.00){\line(1,0){5.00}}
%End Rectangle

\linethickness{0.15mm}
%Rectangle(55.00,10.00)(60.00,15.00)  
\put(55.00,10.00){\line(1,0){5.00}}
\put(55.00,10.00){\line(0,1){5.00}}
\put(60.00,10.00){\line(0,1){5.00}}
\put(55.00,15.00){\line(1,0){5.00}}
%End Rectangle

\linethickness{0.15mm}
%Rectangle(60.00,10.00)(65.00,15.00)  
\put(60.00,10.00){\line(1,0){5.00}}
\put(60.00,10.00){\line(0,1){5.00}}
\put(65.00,10.00){\line(0,1){5.00}}
\put(60.00,15.00){\line(1,0){5.00}}
%End Rectangle

\linethickness{0.15mm}
%Rectangle(65.00,10.00)(70.00,15.00)  
\put(65.00,10.00){\line(1,0){5.00}}
\put(65.00,10.00){\line(0,1){5.00}}
\put(70.00,10.00){\line(0,1){5.00}}
\put(65.00,15.00){\line(1,0){5.00}}
%End Rectangle

\linethickness{0.15mm}
%Rectangle(40.00,15.00)(45.00,20.00)  
\put(40.00,15.00){\line(1,0){5.00}}
\put(40.00,15.00){\line(0,1){5.00}}
\put(45.00,15.00){\line(0,1){5.00}}
\put(40.00,20.00){\line(1,0){5.00}}
%End Rectangle

\linethickness{0.15mm}
%Rectangle(45.00,15.00)(50.00,20.00)  
\put(45.00,15.00){\line(1,0){5.00}}
\put(45.00,15.00){\line(0,1){5.00}}
\put(50.00,15.00){\line(0,1){5.00}}
\put(45.00,20.00){\line(1,0){5.00}}
%End Rectangle

\linethickness{0.15mm}
%Rectangle(50.00,15.00)(55.00,20.00)  
\put(50.00,15.00){\line(1,0){5.00}}
\put(50.00,15.00){\line(0,1){5.00}}
\put(55.00,15.00){\line(0,1){5.00}}
\put(50.00,20.00){\line(1,0){5.00}}
%End Rectangle

\linethickness{0.15mm}
%Rectangle(55.00,15.00)(60.00,20.00)  
\put(55.00,15.00){\line(1,0){5.00}}
\put(55.00,15.00){\line(0,1){5.00}}
\put(60.00,15.00){\line(0,1){5.00}}
\put(55.00,20.00){\line(1,0){5.00}}
%End Rectangle

\linethickness{0.15mm}
%Rectangle(60.00,15.00)(65.00,20.00)  
\put(60.00,15.00){\line(1,0){5.00}}
\put(60.00,15.00){\line(0,1){5.00}}
\put(65.00,15.00){\line(0,1){5.00}}
\put(60.00,20.00){\line(1,0){5.00}}
%End Rectangle

\linethickness{0.15mm}
%Rectangle(40.00,20.00)(45.00,25.00)  
\put(40.00,20.00){\line(1,0){5.00}}
\put(40.00,20.00){\line(0,1){5.00}}
\put(45.00,20.00){\line(0,1){5.00}}
\put(40.00,25.00){\line(1,0){5.00}}
%End Rectangle

\linethickness{0.15mm}
%Rectangle(40.00,25.00)(45.00,30.00)  
\put(40.00,25.00){\line(1,0){5.00}}
\put(40.00,25.00){\line(0,1){5.00}}
\put(45.00,25.00){\line(0,1){5.00}}
\put(40.00,30.00){\line(1,0){5.00}}
%End Rectangle

\linethickness{0.15mm}
%Rectangle(40.00,30.00)(45.00,35.00)  
\put(40.00,30.00){\line(1,0){5.00}}
\put(40.00,30.00){\line(0,1){5.00}}
\put(45.00,30.00){\line(0,1){5.00}}
\put(40.00,35.00){\line(1,0){5.00}}
%End Rectangle

\linethickness{0.15mm}
%Rectangle(45.00,20.00)(50.00,25.00)  
\put(45.00,20.00){\line(1,0){5.00}}
\put(45.00,20.00){\line(0,1){5.00}}
\put(50.00,20.00){\line(0,1){5.00}}
\put(45.00,25.00){\line(1,0){5.00}}
%End Rectangle

\put(35.00,2.50){\makebox(0,0)[cc]{$8$}}

\put(35.00,7.50){\makebox(0,0)[cc]{$6$}}

\put(35.00,12.50){\makebox(0,0)[cc]{$6$}}

\put(35.00,17.50){\makebox(0,0)[cc]{$5$}}

\put(35.00,22.50){\makebox(0,0)[cc]{$2$}}

\put(35.00,27.50){\makebox(0,0)[cc]{$1$}}

\put(35.00,32.50){\makebox(0,0)[cc]{$1$}}

\end{picture} \\

For such a Ferrers graph we consider the chess (or draughts) colouring on it, with the convention that the unit square left below is black. For example the chess Ferrers graph of the partition $\la = (8,6,6,5,2,1,1)$ is given by \\

%\input partition29colour.tex \\
%Created by jPicEdt 1.x
%Standard LaTeX format (emulated lines)
%Fri Jul 08 12:47:38 CEST 2005
\unitlength 1mm
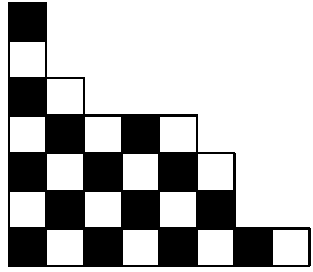
\begin{picture}(80.00,35.00)(0,0)

\linethickness{0.15mm}
%Rectangle(40.00,0.00)(45.00,5.00)  
\put(40.00,0.00){\line(1,0){5.00}}
\put(40.00,0.00){\line(0,1){5.00}}
\put(45.00,0.00){\line(0,1){5.00}}
\put(40.00,5.00){\line(1,0){5.00}}
%End Rectangle

\linethickness{0.15mm}
%Rectangle(45.00,0.00)(50.00,5.00)  
\put(45.00,0.00){\line(1,0){5.00}}
\put(45.00,0.00){\line(0,1){5.00}}
\put(50.00,0.00){\line(0,1){5.00}}
\put(45.00,5.00){\line(1,0){5.00}}
%End Rectangle

\linethickness{0.15mm}
%Rectangle(50.00,0.00)(55.00,5.00)  
\put(50.00,0.00){\line(1,0){5.00}}
\put(50.00,0.00){\line(0,1){5.00}}
\put(55.00,0.00){\line(0,1){5.00}}
\put(50.00,5.00){\line(1,0){5.00}}
%End Rectangle

\linethickness{0.15mm}
%Rectangle(55.00,0.00)(60.00,5.00)  
\put(55.00,0.00){\line(1,0){5.00}}
\put(55.00,0.00){\line(0,1){5.00}}
\put(60.00,0.00){\line(0,1){5.00}}
\put(55.00,5.00){\line(1,0){5.00}}
%End Rectangle

\linethickness{0.15mm}
%Rectangle(60.00,0.00)(65.00,5.00)  
\put(60.00,0.00){\line(1,0){5.00}}
\put(60.00,0.00){\line(0,1){5.00}}
\put(65.00,0.00){\line(0,1){5.00}}
\put(60.00,5.00){\line(1,0){5.00}}
%End Rectangle

\linethickness{0.15mm}
%Rectangle(65.00,0.00)(70.00,5.00)  
\put(65.00,0.00){\line(1,0){5.00}}
\put(65.00,0.00){\line(0,1){5.00}}
\put(70.00,0.00){\line(0,1){5.00}}
\put(65.00,5.00){\line(1,0){5.00}}
%End Rectangle

\linethickness{0.15mm}
%Rectangle(70.00,0.00)(75.00,5.00)  
\put(70.00,0.00){\line(1,0){5.00}}
\put(70.00,0.00){\line(0,1){5.00}}
\put(75.00,0.00){\line(0,1){5.00}}
\put(70.00,5.00){\line(1,0){5.00}}
%End Rectangle

\linethickness{0.15mm}
%Rectangle(75.00,0.00)(80.00,5.00)  
\put(75.00,0.00){\line(1,0){5.00}}
\put(75.00,0.00){\line(0,1){5.00}}
\put(80.00,0.00){\line(0,1){5.00}}
\put(75.00,5.00){\line(1,0){5.00}}
%End Rectangle

\linethickness{0.15mm}
%Rectangle(40.00,5.00)(45.00,10.00)  
\put(40.00,5.00){\line(1,0){5.00}}
\put(40.00,5.00){\line(0,1){5.00}}
\put(45.00,5.00){\line(0,1){5.00}}
\put(40.00,10.00){\line(1,0){5.00}}
%End Rectangle

\linethickness{0.15mm}
%Rectangle(45.00,5.00)(50.00,10.00)  
\put(45.00,5.00){\line(1,0){5.00}}
\put(45.00,5.00){\line(0,1){5.00}}
\put(50.00,5.00){\line(0,1){5.00}}
\put(45.00,10.00){\line(1,0){5.00}}
%End Rectangle

\linethickness{0.15mm}
%Rectangle(50.00,5.00)(55.00,10.00)  
\put(50.00,5.00){\line(1,0){5.00}}
\put(50.00,5.00){\line(0,1){5.00}}
\put(55.00,5.00){\line(0,1){5.00}}
\put(50.00,10.00){\line(1,0){5.00}}
%End Rectangle

\linethickness{0.15mm}
%Rectangle(55.00,5.00)(60.00,10.00)  
\put(55.00,5.00){\line(1,0){5.00}}
\put(55.00,5.00){\line(0,1){5.00}}
\put(60.00,5.00){\line(0,1){5.00}}
\put(55.00,10.00){\line(1,0){5.00}}
%End Rectangle

\linethickness{0.15mm}
%Rectangle(60.00,5.00)(65.00,10.00)  
\put(60.00,5.00){\line(1,0){5.00}}
\put(60.00,5.00){\line(0,1){5.00}}
\put(65.00,5.00){\line(0,1){5.00}}
\put(60.00,10.00){\line(1,0){5.00}}
%End Rectangle

\linethickness{0.15mm}
%Rectangle(65.00,5.00)(70.00,10.00)  
\put(65.00,5.00){\line(1,0){5.00}}
\put(65.00,5.00){\line(0,1){5.00}}
\put(70.00,5.00){\line(0,1){5.00}}
\put(65.00,10.00){\line(1,0){5.00}}
%End Rectangle

\linethickness{0.15mm}
%Rectangle(40.00,10.00)(45.00,15.00)  
\put(40.00,10.00){\line(1,0){5.00}}
\put(40.00,10.00){\line(0,1){5.00}}
\put(45.00,10.00){\line(0,1){5.00}}
\put(40.00,15.00){\line(1,0){5.00}}
%End Rectangle

\linethickness{0.15mm}
%Rectangle(45.00,10.00)(50.00,15.00)  
\put(45.00,10.00){\line(1,0){5.00}}
\put(45.00,10.00){\line(0,1){5.00}}
\put(50.00,10.00){\line(0,1){5.00}}
\put(45.00,15.00){\line(1,0){5.00}}
%End Rectangle

\linethickness{0.15mm}
%Rectangle(50.00,10.00)(55.00,15.00)  
\put(50.00,10.00){\line(1,0){5.00}}
\put(50.00,10.00){\line(0,1){5.00}}
\put(55.00,10.00){\line(0,1){5.00}}
\put(50.00,15.00){\line(1,0){5.00}}
%End Rectangle

\linethickness{0.15mm}
%Rectangle(55.00,10.00)(60.00,15.00)  
\put(55.00,10.00){\line(1,0){5.00}}
\put(55.00,10.00){\line(0,1){5.00}}
\put(60.00,10.00){\line(0,1){5.00}}
\put(55.00,15.00){\line(1,0){5.00}}
%End Rectangle

\linethickness{0.15mm}
%Rectangle(60.00,10.00)(65.00,15.00)  
\put(60.00,10.00){\line(1,0){5.00}}
\put(60.00,10.00){\line(0,1){5.00}}
\put(65.00,10.00){\line(0,1){5.00}}
\put(60.00,15.00){\line(1,0){5.00}}
%End Rectangle

\linethickness{0.15mm}
%Rectangle(65.00,10.00)(70.00,15.00)  
\put(65.00,10.00){\line(1,0){5.00}}
\put(65.00,10.00){\line(0,1){5.00}}
\put(70.00,10.00){\line(0,1){5.00}}
\put(65.00,15.00){\line(1,0){5.00}}
%End Rectangle

\linethickness{0.15mm}
%Rectangle(40.00,15.00)(45.00,20.00)  
\put(40.00,15.00){\line(1,0){5.00}}
\put(40.00,15.00){\line(0,1){5.00}}
\put(45.00,15.00){\line(0,1){5.00}}
\put(40.00,20.00){\line(1,0){5.00}}
%End Rectangle

\linethickness{0.15mm}
%Rectangle(45.00,15.00)(50.00,20.00)  
\put(45.00,15.00){\line(1,0){5.00}}
\put(45.00,15.00){\line(0,1){5.00}}
\put(50.00,15.00){\line(0,1){5.00}}
\put(45.00,20.00){\line(1,0){5.00}}
%End Rectangle

\linethickness{0.15mm}
%Rectangle(50.00,15.00)(55.00,20.00)  
\put(50.00,15.00){\line(1,0){5.00}}
\put(50.00,15.00){\line(0,1){5.00}}
\put(55.00,15.00){\line(0,1){5.00}}
\put(50.00,20.00){\line(1,0){5.00}}
%End Rectangle

\linethickness{0.15mm}
%Rectangle(55.00,15.00)(60.00,20.00)  
\put(55.00,15.00){\line(1,0){5.00}}
\put(55.00,15.00){\line(0,1){5.00}}
\put(60.00,15.00){\line(0,1){5.00}}
\put(55.00,20.00){\line(1,0){5.00}}
%End Rectangle

\linethickness{0.15mm}
%Rectangle(60.00,15.00)(65.00,20.00)  
\put(60.00,15.00){\line(1,0){5.00}}
\put(60.00,15.00){\line(0,1){5.00}}
\put(65.00,15.00){\line(0,1){5.00}}
\put(60.00,20.00){\line(1,0){5.00}}
%End Rectangle

\linethickness{0.15mm}
%Rectangle(40.00,20.00)(45.00,25.00)  
\put(40.00,20.00){\line(1,0){5.00}}
\put(40.00,20.00){\line(0,1){5.00}}
\put(45.00,20.00){\line(0,1){5.00}}
\put(40.00,25.00){\line(1,0){5.00}}
%End Rectangle

\linethickness{0.15mm}
%Rectangle(40.00,25.00)(45.00,30.00)  
\put(40.00,25.00){\line(1,0){5.00}}
\put(40.00,25.00){\line(0,1){5.00}}
\put(45.00,25.00){\line(0,1){5.00}}
\put(40.00,30.00){\line(1,0){5.00}}
%End Rectangle

\linethickness{0.15mm}
%Rectangle(40.00,30.00)(45.00,35.00)  
\put(40.00,30.00){\line(1,0){5.00}}
\put(40.00,30.00){\line(0,1){5.00}}
\put(45.00,30.00){\line(0,1){5.00}}
\put(40.00,35.00){\line(1,0){5.00}}
%End Rectangle

\linethickness{0.15mm}
%Rectangle(40.00,0.00)(45.00,5.00)  blacken
\put(40.00,0.00){\rule{5.00\unitlength}{5.00\unitlength}}
%End Rectangle

\linethickness{0.15mm}
%Rectangle(45.00,5.00)(50.00,10.00)  blacken
\put(45.00,5.00){\rule{5.00\unitlength}{5.00\unitlength}}
%End Rectangle

\linethickness{0.15mm}
%Rectangle(50.00,10.00)(55.00,15.00)  blacken
\put(50.00,10.00){\rule{5.00\unitlength}{5.00\unitlength}}
%End Rectangle

\linethickness{0.15mm}
%Rectangle(55.00,15.00)(60.00,20.00)  blacken
\put(55.00,15.00){\rule{5.00\unitlength}{5.00\unitlength}}
%End Rectangle

\linethickness{0.15mm}
%Rectangle(50.00,0.00)(55.00,5.00)  blacken
\put(50.00,0.00){\rule{5.00\unitlength}{5.00\unitlength}}
%End Rectangle

\linethickness{0.15mm}
%Rectangle(55.00,5.00)(60.00,10.00)  blacken
\put(55.00,5.00){\rule{5.00\unitlength}{5.00\unitlength}}
%End Rectangle

\linethickness{0.15mm}
%Rectangle(60.00,10.00)(65.00,15.00)  blacken
\put(60.00,10.00){\rule{5.00\unitlength}{5.00\unitlength}}
%End Rectangle

\linethickness{0.15mm}
%Rectangle(60.00,0.00)(65.00,5.00)  blacken
\put(60.00,0.00){\rule{5.00\unitlength}{5.00\unitlength}}
%End Rectangle

\linethickness{0.15mm}
%Rectangle(65.00,5.00)(70.00,10.00)  blacken
\put(65.00,5.00){\rule{5.00\unitlength}{5.00\unitlength}}
%End Rectangle

\linethickness{0.15mm}
%Rectangle(70.00,0.00)(75.00,5.00)  blacken
\put(70.00,0.00){\rule{5.00\unitlength}{5.00\unitlength}}
%End Rectangle

\linethickness{0.15mm}
%Rectangle(40.00,10.00)(45.00,15.00)  blacken
\put(40.00,10.00){\rule{5.00\unitlength}{5.00\unitlength}}
%End Rectangle

\linethickness{0.15mm}
%Rectangle(45.00,15.00)(50.00,20.00)  blacken
\put(45.00,15.00){\rule{5.00\unitlength}{5.00\unitlength}}
%End Rectangle

\linethickness{0.15mm}
%Rectangle(40.00,20.00)(45.00,25.00)  blacken
\put(40.00,20.00){\rule{5.00\unitlength}{5.00\unitlength}}
%End Rectangle

\linethickness{0.15mm}
%Rectangle(40.00,30.00)(45.00,35.00)  blacken
\put(40.00,30.00){\rule{5.00\unitlength}{5.00\unitlength}}
%End Rectangle

\linethickness{0.15mm}
%Rectangle(45.00,20.00)(50.00,25.00)  
\put(45.00,20.00){\line(1,0){5.00}}
\put(45.00,20.00){\line(0,1){5.00}}
\put(50.00,20.00){\line(0,1){5.00}}
\put(45.00,25.00){\line(1,0){5.00}}
%End Rectangle

\linethickness{0.15mm}
%Rectangle(40.00,25.00)(45.00,30.00)  
\put(40.00,25.00){\line(1,0){5.00}}
\put(40.00,25.00){\line(0,1){5.00}}
\put(45.00,25.00){\line(0,1){5.00}}
\put(40.00,30.00){\line(1,0){5.00}}
%End Rectangle

\linethickness{0.15mm}
%Rectangle(45.00,20.00)(50.00,25.00)  
\put(45.00,20.00){\line(1,0){5.00}}
\put(45.00,20.00){\line(0,1){5.00}}
\put(50.00,20.00){\line(0,1){5.00}}
\put(45.00,25.00){\line(1,0){5.00}}
%End Rectangle

\linethickness{0.15mm}
%Rectangle(50.00,15.00)(55.00,20.00)  
\put(50.00,15.00){\line(1,0){5.00}}
\put(50.00,15.00){\line(0,1){5.00}}
\put(55.00,15.00){\line(0,1){5.00}}
\put(50.00,20.00){\line(1,0){5.00}}
%End Rectangle

\linethickness{0.15mm}
%Rectangle(60.00,15.00)(65.00,20.00)  
\put(60.00,15.00){\line(1,0){5.00}}
\put(60.00,15.00){\line(0,1){5.00}}
\put(65.00,15.00){\line(0,1){5.00}}
\put(60.00,20.00){\line(1,0){5.00}}
%End Rectangle

\linethickness{0.15mm}
%Rectangle(65.00,10.00)(70.00,15.00)  
\put(65.00,10.00){\line(1,0){5.00}}
\put(65.00,10.00){\line(0,1){5.00}}
\put(70.00,10.00){\line(0,1){5.00}}
\put(65.00,15.00){\line(1,0){5.00}}
%End Rectangle

\linethickness{0.15mm}
%Rectangle(75.00,0.00)(80.00,5.00)  
\put(75.00,0.00){\line(1,0){5.00}}
\put(75.00,0.00){\line(0,1){5.00}}
\put(80.00,0.00){\line(0,1){5.00}}
\put(75.00,5.00){\line(1,0){5.00}}
%End Rectangle

\end{picture} \\

For a partition $\la$ we write $b(\la)$ (resp. $w(\la)$) for the number of black (resp. white) squares in its chess Ferrers graph. Our main result is
\begin{theoremA} 
Let $(b,w) \in \NN^{2}$. Then there exists a partition $\la$ such that \linebreak $(b(\la),w(\la)) = (b,w)$ if and only if 
\begin{equation} \label{ineq}
(b-w)^{2} \leq b
\end{equation}
Furthermore the same statement holds if we restrict ourselves to partitions in distinct parts.
\end{theoremA}
If $b \neq 0$ then \eqref{ineq} may be written as
\[
\left(1-\frac{w}{b}\right)^{2} \leq \frac{1}{b}
\]
%which we might call a broken chessboard inequality. 
which measures how close the ratio $w/b$ is to $1$.
As a byproduct of the proof of Theorem A presented in this note, the appearing $(b,w) \in \NN^{2}$ are discribed in an explicit way:
\begin{theoremB} 
Let $(b,w) \in \NN^{2}$. Then there exists a partition $\la$ such that \linebreak $(b(\la),w(\la)) = (b,w)$ if and only if there exist positive integers $k,l \in \NN$ such that either
\[
(b,w) = \left( (k+1)^{2} + l ,k(k+1) + l \right) \text{ or } (b,w) = \left( k^{2} + l, k(k+1) + l \right)
\]
\end{theoremB}
Let us indicate intuitively how we prove Theorem A. To any chess Ferrers graph we associate another graph by  
\begin{enumerate}
\item shifting the first row one place to the right, the second row two places to the right, etc. and afterwards
\item
if necessary filling the ``holes" by applying gravity.
\end{enumerate}
For example for the partition $\la = (8,6,6,5,2,1,1)$ we find \\

%\input procedure29.tex \\
%Created by jPicEdt 1.x
%Standard LaTeX format (emulated lines)
%Wed Jul 06 22:17:06 CEST 2005
\unitlength 0.9mm
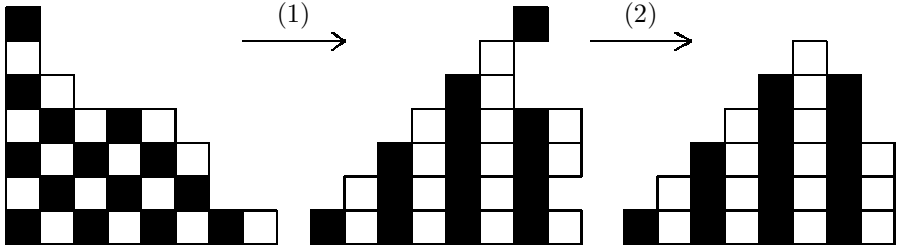
\begin{picture}(126.25,35.00)(0,0)

\linethickness{0.15mm}
%Rectangle(-5.00,0.00)(0.00,5.00)  
\put(-5.00,0.00){\line(1,0){5.00}}
\put(-5.00,0.00){\line(0,1){5.00}}
\put(0.00,0.00){\line(0,1){5.00}}
\put(-5.00,5.00){\line(1,0){5.00}}
%End Rectangle

\linethickness{0.15mm}
%Rectangle(0.00,0.00)(5.00,5.00)  
\put(0.00,0.00){\line(1,0){5.00}}
\put(0.00,0.00){\line(0,1){5.00}}
\put(5.00,0.00){\line(0,1){5.00}}
\put(0.00,5.00){\line(1,0){5.00}}
%End Rectangle

\linethickness{0.15mm}
%Rectangle(5.00,0.00)(10.00,5.00)  
\put(5.00,0.00){\line(1,0){5.00}}
\put(5.00,0.00){\line(0,1){5.00}}
\put(10.00,0.00){\line(0,1){5.00}}
\put(5.00,5.00){\line(1,0){5.00}}
%End Rectangle

\linethickness{0.15mm}
%Rectangle(10.00,0.00)(15.00,5.00)  
\put(10.00,0.00){\line(1,0){5.00}}
\put(10.00,0.00){\line(0,1){5.00}}
\put(15.00,0.00){\line(0,1){5.00}}
\put(10.00,5.00){\line(1,0){5.00}}
%End Rectangle

\linethickness{0.15mm}
%Rectangle(15.00,0.00)(20.00,5.00)  
\put(15.00,0.00){\line(1,0){5.00}}
\put(15.00,0.00){\line(0,1){5.00}}
\put(20.00,0.00){\line(0,1){5.00}}
\put(15.00,5.00){\line(1,0){5.00}}
%End Rectangle

\linethickness{0.15mm}
%Rectangle(20.00,0.00)(25.00,5.00)  
\put(20.00,0.00){\line(1,0){5.00}}
\put(20.00,0.00){\line(0,1){5.00}}
\put(25.00,0.00){\line(0,1){5.00}}
\put(20.00,5.00){\line(1,0){5.00}}
%End Rectangle

\linethickness{0.15mm}
%Rectangle(25.00,0.00)(30.00,5.00)  
\put(25.00,0.00){\line(1,0){5.00}}
\put(25.00,0.00){\line(0,1){5.00}}
\put(30.00,0.00){\line(0,1){5.00}}
\put(25.00,5.00){\line(1,0){5.00}}
%End Rectangle

\linethickness{0.15mm}
%Rectangle(30.00,0.00)(35.00,5.00)  
\put(30.00,0.00){\line(1,0){5.00}}
\put(30.00,0.00){\line(0,1){5.00}}
\put(35.00,0.00){\line(0,1){5.00}}
\put(30.00,5.00){\line(1,0){5.00}}
%End Rectangle

\linethickness{0.15mm}
%Rectangle(-5.00,5.00)(0.00,10.00)  
\put(-5.00,5.00){\line(1,0){5.00}}
\put(-5.00,5.00){\line(0,1){5.00}}
\put(0.00,5.00){\line(0,1){5.00}}
\put(-5.00,10.00){\line(1,0){5.00}}
%End Rectangle

\linethickness{0.15mm}
%Rectangle(0.00,5.00)(5.00,10.00)  
\put(0.00,5.00){\line(1,0){5.00}}
\put(0.00,5.00){\line(0,1){5.00}}
\put(5.00,5.00){\line(0,1){5.00}}
\put(0.00,10.00){\line(1,0){5.00}}
%End Rectangle

\linethickness{0.15mm}
%Rectangle(5.00,5.00)(10.00,10.00)  
\put(5.00,5.00){\line(1,0){5.00}}
\put(5.00,5.00){\line(0,1){5.00}}
\put(10.00,5.00){\line(0,1){5.00}}
\put(5.00,10.00){\line(1,0){5.00}}
%End Rectangle

\linethickness{0.15mm}
%Rectangle(10.00,5.00)(15.00,10.00)  
\put(10.00,5.00){\line(1,0){5.00}}
\put(10.00,5.00){\line(0,1){5.00}}
\put(15.00,5.00){\line(0,1){5.00}}
\put(10.00,10.00){\line(1,0){5.00}}
%End Rectangle

\linethickness{0.15mm}
%Rectangle(15.00,5.00)(20.00,10.00)  
\put(15.00,5.00){\line(1,0){5.00}}
\put(15.00,5.00){\line(0,1){5.00}}
\put(20.00,5.00){\line(0,1){5.00}}
\put(15.00,10.00){\line(1,0){5.00}}
%End Rectangle

\linethickness{0.15mm}
%Rectangle(20.00,5.00)(25.00,10.00)  
\put(20.00,5.00){\line(1,0){5.00}}
\put(20.00,5.00){\line(0,1){5.00}}
\put(25.00,5.00){\line(0,1){5.00}}
\put(20.00,10.00){\line(1,0){5.00}}
%End Rectangle

\linethickness{0.15mm}
%Rectangle(-5.00,10.00)(0.00,15.00)  
\put(-5.00,10.00){\line(1,0){5.00}}
\put(-5.00,10.00){\line(0,1){5.00}}
\put(0.00,10.00){\line(0,1){5.00}}
\put(-5.00,15.00){\line(1,0){5.00}}
%End Rectangle

\linethickness{0.15mm}
%Rectangle(0.00,10.00)(5.00,15.00)  
\put(0.00,10.00){\line(1,0){5.00}}
\put(0.00,10.00){\line(0,1){5.00}}
\put(5.00,10.00){\line(0,1){5.00}}
\put(0.00,15.00){\line(1,0){5.00}}
%End Rectangle

\linethickness{0.15mm}
%Rectangle(5.00,10.00)(10.00,15.00)  
\put(5.00,10.00){\line(1,0){5.00}}
\put(5.00,10.00){\line(0,1){5.00}}
\put(10.00,10.00){\line(0,1){5.00}}
\put(5.00,15.00){\line(1,0){5.00}}
%End Rectangle

\linethickness{0.15mm}
%Rectangle(10.00,10.00)(15.00,15.00)  
\put(10.00,10.00){\line(1,0){5.00}}
\put(10.00,10.00){\line(0,1){5.00}}
\put(15.00,10.00){\line(0,1){5.00}}
\put(10.00,15.00){\line(1,0){5.00}}
%End Rectangle

\linethickness{0.15mm}
%Rectangle(15.00,10.00)(20.00,15.00)  
\put(15.00,10.00){\line(1,0){5.00}}
\put(15.00,10.00){\line(0,1){5.00}}
\put(20.00,10.00){\line(0,1){5.00}}
\put(15.00,15.00){\line(1,0){5.00}}
%End Rectangle

\linethickness{0.15mm}
%Rectangle(20.00,10.00)(25.00,15.00)  
\put(20.00,10.00){\line(1,0){5.00}}
\put(20.00,10.00){\line(0,1){5.00}}
\put(25.00,10.00){\line(0,1){5.00}}
\put(20.00,15.00){\line(1,0){5.00}}
%End Rectangle

\linethickness{0.15mm}
%Rectangle(-5.00,15.00)(0.00,20.00)  
\put(-5.00,15.00){\line(1,0){5.00}}
\put(-5.00,15.00){\line(0,1){5.00}}
\put(0.00,15.00){\line(0,1){5.00}}
\put(-5.00,20.00){\line(1,0){5.00}}
%End Rectangle

\linethickness{0.15mm}
%Rectangle(0.00,15.00)(5.00,20.00)  
\put(0.00,15.00){\line(1,0){5.00}}
\put(0.00,15.00){\line(0,1){5.00}}
\put(5.00,15.00){\line(0,1){5.00}}
\put(0.00,20.00){\line(1,0){5.00}}
%End Rectangle

\linethickness{0.15mm}
%Rectangle(5.00,15.00)(10.00,20.00)  
\put(5.00,15.00){\line(1,0){5.00}}
\put(5.00,15.00){\line(0,1){5.00}}
\put(10.00,15.00){\line(0,1){5.00}}
\put(5.00,20.00){\line(1,0){5.00}}
%End Rectangle

\linethickness{0.15mm}
%Rectangle(10.00,15.00)(15.00,20.00)  
\put(10.00,15.00){\line(1,0){5.00}}
\put(10.00,15.00){\line(0,1){5.00}}
\put(15.00,15.00){\line(0,1){5.00}}
\put(10.00,20.00){\line(1,0){5.00}}
%End Rectangle

\linethickness{0.15mm}
%Rectangle(15.00,15.00)(20.00,20.00)  
\put(15.00,15.00){\line(1,0){5.00}}
\put(15.00,15.00){\line(0,1){5.00}}
\put(20.00,15.00){\line(0,1){5.00}}
\put(15.00,20.00){\line(1,0){5.00}}
%End Rectangle

\linethickness{0.15mm}
%Rectangle(-5.00,20.00)(0.00,25.00)  
\put(-5.00,20.00){\line(1,0){5.00}}
\put(-5.00,20.00){\line(0,1){5.00}}
\put(0.00,20.00){\line(0,1){5.00}}
\put(-5.00,25.00){\line(1,0){5.00}}
%End Rectangle

\linethickness{0.15mm}
%Rectangle(-5.00,25.00)(0.00,30.00)  
\put(-5.00,25.00){\line(1,0){5.00}}
\put(-5.00,25.00){\line(0,1){5.00}}
\put(0.00,25.00){\line(0,1){5.00}}
\put(-5.00,30.00){\line(1,0){5.00}}
%End Rectangle

\linethickness{0.15mm}
%Rectangle(-5.00,30.00)(0.00,35.00)  
\put(-5.00,30.00){\line(1,0){5.00}}
\put(-5.00,30.00){\line(0,1){5.00}}
\put(0.00,30.00){\line(0,1){5.00}}
\put(-5.00,35.00){\line(1,0){5.00}}
%End Rectangle

\linethickness{0.15mm}
%Rectangle(-5.00,0.00)(0.00,5.00)  blacken
\put(-5.00,0.00){\rule{5.00\unitlength}{5.00\unitlength}}
%End Rectangle

\linethickness{0.15mm}
%Rectangle(0.00,5.00)(5.00,10.00)  blacken
\put(0.00,5.00){\rule{5.00\unitlength}{5.00\unitlength}}
%End Rectangle

\linethickness{0.15mm}
%Rectangle(5.00,10.00)(10.00,15.00)  blacken
\put(5.00,10.00){\rule{5.00\unitlength}{5.00\unitlength}}
%End Rectangle

\linethickness{0.15mm}
%Rectangle(10.00,15.00)(15.00,20.00)  blacken
\put(10.00,15.00){\rule{5.00\unitlength}{5.00\unitlength}}
%End Rectangle

\linethickness{0.15mm}
%Rectangle(5.00,0.00)(10.00,5.00)  blacken
\put(5.00,0.00){\rule{5.00\unitlength}{5.00\unitlength}}
%End Rectangle

\linethickness{0.15mm}
%Rectangle(10.00,5.00)(15.00,10.00)  blacken
\put(10.00,5.00){\rule{5.00\unitlength}{5.00\unitlength}}
%End Rectangle

\linethickness{0.15mm}
%Rectangle(15.00,10.00)(20.00,15.00)  blacken
\put(15.00,10.00){\rule{5.00\unitlength}{5.00\unitlength}}
%End Rectangle

\linethickness{0.15mm}
%Rectangle(15.00,0.00)(20.00,5.00)  blacken
\put(15.00,0.00){\rule{5.00\unitlength}{5.00\unitlength}}
%End Rectangle

\linethickness{0.15mm}
%Rectangle(20.00,5.00)(25.00,10.00)  blacken
\put(20.00,5.00){\rule{5.00\unitlength}{5.00\unitlength}}
%End Rectangle

\linethickness{0.15mm}
%Rectangle(25.00,0.00)(30.00,5.00)  blacken
\put(25.00,0.00){\rule{5.00\unitlength}{5.00\unitlength}}
%End Rectangle

\linethickness{0.15mm}
%Rectangle(-5.00,10.00)(0.00,15.00)  blacken
\put(-5.00,10.00){\rule{5.00\unitlength}{5.00\unitlength}}
%End Rectangle

\linethickness{0.15mm}
%Rectangle(0.00,15.00)(5.00,20.00)  blacken
\put(0.00,15.00){\rule{5.00\unitlength}{5.00\unitlength}}
%End Rectangle

\linethickness{0.15mm}
%Rectangle(-5.00,20.00)(0.00,25.00)  blacken
\put(-5.00,20.00){\rule{5.00\unitlength}{5.00\unitlength}}
%End Rectangle

\linethickness{0.15mm}
%Rectangle(-5.00,30.00)(0.00,35.00)  blacken
\put(-5.00,30.00){\rule{5.00\unitlength}{5.00\unitlength}}
%End Rectangle

\linethickness{0.15mm}
%Rectangle(0.00,20.00)(5.00,25.00)  
\put(0.00,20.00){\line(1,0){5.00}}
\put(0.00,20.00){\line(0,1){5.00}}
\put(5.00,20.00){\line(0,1){5.00}}
\put(0.00,25.00){\line(1,0){5.00}}
%End Rectangle

\linethickness{0.15mm}
%Rectangle(40.00,0.00)(45.00,5.00)  
\put(40.00,0.00){\line(1,0){5.00}}
\put(40.00,0.00){\line(0,1){5.00}}
\put(45.00,0.00){\line(0,1){5.00}}
\put(40.00,5.00){\line(1,0){5.00}}
%End Rectangle

\linethickness{0.15mm}
%Rectangle(45.00,0.00)(50.00,5.00)  
\put(45.00,0.00){\line(1,0){5.00}}
\put(45.00,0.00){\line(0,1){5.00}}
\put(50.00,0.00){\line(0,1){5.00}}
\put(45.00,5.00){\line(1,0){5.00}}
%End Rectangle

\linethickness{0.15mm}
%Rectangle(50.00,0.00)(55.00,5.00)  
\put(50.00,0.00){\line(1,0){5.00}}
\put(50.00,0.00){\line(0,1){5.00}}
\put(55.00,0.00){\line(0,1){5.00}}
\put(50.00,5.00){\line(1,0){5.00}}
%End Rectangle

\linethickness{0.15mm}
%Rectangle(55.00,0.00)(60.00,5.00)  
\put(55.00,0.00){\line(1,0){5.00}}
\put(55.00,0.00){\line(0,1){5.00}}
\put(60.00,0.00){\line(0,1){5.00}}
\put(55.00,5.00){\line(1,0){5.00}}
%End Rectangle

\linethickness{0.15mm}
%Rectangle(60.00,0.00)(65.00,5.00)  
\put(60.00,0.00){\line(1,0){5.00}}
\put(60.00,0.00){\line(0,1){5.00}}
\put(65.00,0.00){\line(0,1){5.00}}
\put(60.00,5.00){\line(1,0){5.00}}
%End Rectangle

\linethickness{0.15mm}
%Rectangle(65.00,0.00)(70.00,5.00)  
\put(65.00,0.00){\line(1,0){5.00}}
\put(65.00,0.00){\line(0,1){5.00}}
\put(70.00,0.00){\line(0,1){5.00}}
\put(65.00,5.00){\line(1,0){5.00}}
%End Rectangle

\linethickness{0.15mm}
%Rectangle(70.00,0.00)(75.00,5.00)  
\put(70.00,0.00){\line(1,0){5.00}}
\put(70.00,0.00){\line(0,1){5.00}}
\put(75.00,0.00){\line(0,1){5.00}}
\put(70.00,5.00){\line(1,0){5.00}}
%End Rectangle

\linethickness{0.15mm}
%Rectangle(75.00,0.00)(80.00,5.00)  
\put(75.00,0.00){\line(1,0){5.00}}
\put(75.00,0.00){\line(0,1){5.00}}
\put(80.00,0.00){\line(0,1){5.00}}
\put(75.00,5.00){\line(1,0){5.00}}
%End Rectangle

\linethickness{0.15mm}
%Rectangle(45.00,5.00)(50.00,10.00)  
\put(45.00,5.00){\line(1,0){5.00}}
\put(45.00,5.00){\line(0,1){5.00}}
\put(50.00,5.00){\line(0,1){5.00}}
\put(45.00,10.00){\line(1,0){5.00}}
%End Rectangle

\linethickness{0.15mm}
%Rectangle(50.00,5.00)(55.00,10.00)  
\put(50.00,5.00){\line(1,0){5.00}}
\put(50.00,5.00){\line(0,1){5.00}}
\put(55.00,5.00){\line(0,1){5.00}}
\put(50.00,10.00){\line(1,0){5.00}}
%End Rectangle

\linethickness{0.15mm}
%Rectangle(55.00,5.00)(60.00,10.00)  
\put(55.00,5.00){\line(1,0){5.00}}
\put(55.00,5.00){\line(0,1){5.00}}
\put(60.00,5.00){\line(0,1){5.00}}
\put(55.00,10.00){\line(1,0){5.00}}
%End Rectangle

\linethickness{0.15mm}
%Rectangle(60.00,5.00)(65.00,10.00)  
\put(60.00,5.00){\line(1,0){5.00}}
\put(60.00,5.00){\line(0,1){5.00}}
\put(65.00,5.00){\line(0,1){5.00}}
\put(60.00,10.00){\line(1,0){5.00}}
%End Rectangle

\linethickness{0.15mm}
%Rectangle(65.00,5.00)(70.00,10.00)  
\put(65.00,5.00){\line(1,0){5.00}}
\put(65.00,5.00){\line(0,1){5.00}}
\put(70.00,5.00){\line(0,1){5.00}}
\put(65.00,10.00){\line(1,0){5.00}}
%End Rectangle

\linethickness{0.15mm}
%Rectangle(70.00,5.00)(75.00,10.00)  
\put(70.00,5.00){\line(1,0){5.00}}
\put(70.00,5.00){\line(0,1){5.00}}
\put(75.00,5.00){\line(0,1){5.00}}
\put(70.00,10.00){\line(1,0){5.00}}
%End Rectangle

\linethickness{0.15mm}
%Rectangle(50.00,10.00)(55.00,15.00)  
\put(50.00,10.00){\line(1,0){5.00}}
\put(50.00,10.00){\line(0,1){5.00}}
\put(55.00,10.00){\line(0,1){5.00}}
\put(50.00,15.00){\line(1,0){5.00}}
%End Rectangle

\linethickness{0.15mm}
%Rectangle(55.00,10.00)(60.00,15.00)  
\put(55.00,10.00){\line(1,0){5.00}}
\put(55.00,10.00){\line(0,1){5.00}}
\put(60.00,10.00){\line(0,1){5.00}}
\put(55.00,15.00){\line(1,0){5.00}}
%End Rectangle

\linethickness{0.15mm}
%Rectangle(60.00,10.00)(65.00,15.00)  
\put(60.00,10.00){\line(1,0){5.00}}
\put(60.00,10.00){\line(0,1){5.00}}
\put(65.00,10.00){\line(0,1){5.00}}
\put(60.00,15.00){\line(1,0){5.00}}
%End Rectangle

\linethickness{0.15mm}
%Rectangle(65.00,10.00)(70.00,15.00)  
\put(65.00,10.00){\line(1,0){5.00}}
\put(65.00,10.00){\line(0,1){5.00}}
\put(70.00,10.00){\line(0,1){5.00}}
\put(65.00,15.00){\line(1,0){5.00}}
%End Rectangle

\linethickness{0.15mm}
%Rectangle(70.00,10.00)(75.00,15.00)  
\put(70.00,10.00){\line(1,0){5.00}}
\put(70.00,10.00){\line(0,1){5.00}}
\put(75.00,10.00){\line(0,1){5.00}}
\put(70.00,15.00){\line(1,0){5.00}}
%End Rectangle

\linethickness{0.15mm}
%Rectangle(75.00,10.00)(80.00,15.00)  
\put(75.00,10.00){\line(1,0){5.00}}
\put(75.00,10.00){\line(0,1){5.00}}
\put(80.00,10.00){\line(0,1){5.00}}
\put(75.00,15.00){\line(1,0){5.00}}
%End Rectangle

\linethickness{0.15mm}
%Rectangle(55.00,15.00)(60.00,20.00)  
\put(55.00,15.00){\line(1,0){5.00}}
\put(55.00,15.00){\line(0,1){5.00}}
\put(60.00,15.00){\line(0,1){5.00}}
\put(55.00,20.00){\line(1,0){5.00}}
%End Rectangle

\linethickness{0.15mm}
%Rectangle(60.00,15.00)(65.00,20.00)  
\put(60.00,15.00){\line(1,0){5.00}}
\put(60.00,15.00){\line(0,1){5.00}}
\put(65.00,15.00){\line(0,1){5.00}}
\put(60.00,20.00){\line(1,0){5.00}}
%End Rectangle

\linethickness{0.15mm}
%Rectangle(65.00,15.00)(70.00,20.00)  
\put(65.00,15.00){\line(1,0){5.00}}
\put(65.00,15.00){\line(0,1){5.00}}
\put(70.00,15.00){\line(0,1){5.00}}
\put(65.00,20.00){\line(1,0){5.00}}
%End Rectangle

\linethickness{0.15mm}
%Rectangle(70.00,15.00)(75.00,20.00)  
\put(70.00,15.00){\line(1,0){5.00}}
\put(70.00,15.00){\line(0,1){5.00}}
\put(75.00,15.00){\line(0,1){5.00}}
\put(70.00,20.00){\line(1,0){5.00}}
%End Rectangle

\linethickness{0.15mm}
%Rectangle(75.00,15.00)(80.00,20.00)  
\put(75.00,15.00){\line(1,0){5.00}}
\put(75.00,15.00){\line(0,1){5.00}}
\put(80.00,15.00){\line(0,1){5.00}}
\put(75.00,20.00){\line(1,0){5.00}}
%End Rectangle

\linethickness{0.15mm}
%Rectangle(60.00,20.00)(65.00,25.00)  
\put(60.00,20.00){\line(1,0){5.00}}
\put(60.00,20.00){\line(0,1){5.00}}
\put(65.00,20.00){\line(0,1){5.00}}
\put(60.00,25.00){\line(1,0){5.00}}
%End Rectangle

\linethickness{0.15mm}
%Rectangle(65.00,20.00)(70.00,25.00)  
\put(65.00,20.00){\line(1,0){5.00}}
\put(65.00,20.00){\line(0,1){5.00}}
\put(70.00,20.00){\line(0,1){5.00}}
\put(65.00,25.00){\line(1,0){5.00}}
%End Rectangle

\linethickness{0.15mm}
%Rectangle(65.00,25.00)(70.00,30.00)  
\put(65.00,25.00){\line(1,0){5.00}}
\put(65.00,25.00){\line(0,1){5.00}}
\put(70.00,25.00){\line(0,1){5.00}}
\put(65.00,30.00){\line(1,0){5.00}}
%End Rectangle

\linethickness{0.15mm}
%Rectangle(70.00,30.00)(75.00,35.00)  
\put(70.00,30.00){\line(1,0){5.00}}
\put(70.00,30.00){\line(0,1){5.00}}
\put(75.00,30.00){\line(0,1){5.00}}
\put(70.00,35.00){\line(1,0){5.00}}
%End Rectangle

\linethickness{0.15mm}
%Rectangle(40.00,0.00)(45.00,5.00)  blacken
\put(40.00,0.00){\rule{5.00\unitlength}{5.00\unitlength}}
%End Rectangle

\linethickness{0.15mm}
%Rectangle(50.00,10.00)(55.00,15.00)  blacken
\put(50.00,10.00){\rule{5.00\unitlength}{5.00\unitlength}}
%End Rectangle

\linethickness{0.15mm}
%Rectangle(50.00,5.00)(55.00,10.00)  blacken
\put(50.00,5.00){\rule{5.00\unitlength}{5.00\unitlength}}
%End Rectangle

\linethickness{0.15mm}
%Rectangle(50.00,0.00)(55.00,5.00)  blacken
\put(50.00,0.00){\rule{5.00\unitlength}{5.00\unitlength}}
%End Rectangle

\linethickness{0.15mm}
%Rectangle(60.00,0.00)(65.00,25.00)  blacken
\put(60.00,0.00){\rule{5.00\unitlength}{25.00\unitlength}}
%End Rectangle

\linethickness{0.15mm}
%Rectangle(70.00,35.00)(75.00,35.00)  blacken
\put(70.00,35.00){\rule{5.00\unitlength}{0.00\unitlength}}
%End Rectangle

\linethickness{0.15mm}
%Rectangle(70.00,30.00)(75.00,35.00)  blacken
\put(70.00,30.00){\rule{5.00\unitlength}{5.00\unitlength}}
%End Rectangle

\linethickness{0.15mm}
%Rectangle(70.00,0.00)(75.00,20.00)  blacken
\put(70.00,0.00){\rule{5.00\unitlength}{20.00\unitlength}}
%End Rectangle

\linethickness{0.15mm}
%Rectangle(86.25,0.00)(91.25,5.00)  
\put(86.25,0.00){\line(1,0){5.00}}
\put(86.25,0.00){\line(0,1){5.00}}
\put(91.25,0.00){\line(0,1){5.00}}
\put(86.25,5.00){\line(1,0){5.00}}
%End Rectangle

\linethickness{0.15mm}
%Rectangle(91.25,0.00)(96.25,5.00)  
\put(91.25,0.00){\line(1,0){5.00}}
\put(91.25,0.00){\line(0,1){5.00}}
\put(96.25,0.00){\line(0,1){5.00}}
\put(91.25,5.00){\line(1,0){5.00}}
%End Rectangle

\linethickness{0.15mm}
%Rectangle(96.25,0.00)(101.25,5.00)  
\put(96.25,0.00){\line(1,0){5.00}}
\put(96.25,0.00){\line(0,1){5.00}}
\put(101.25,0.00){\line(0,1){5.00}}
\put(96.25,5.00){\line(1,0){5.00}}
%End Rectangle

\linethickness{0.15mm}
%Rectangle(101.25,0.00)(106.25,5.00)  
\put(101.25,0.00){\line(1,0){5.00}}
\put(101.25,0.00){\line(0,1){5.00}}
\put(106.25,0.00){\line(0,1){5.00}}
\put(101.25,5.00){\line(1,0){5.00}}
%End Rectangle

\linethickness{0.15mm}
%Rectangle(106.25,0.00)(111.25,5.00)  
\put(106.25,0.00){\line(1,0){5.00}}
\put(106.25,0.00){\line(0,1){5.00}}
\put(111.25,0.00){\line(0,1){5.00}}
\put(106.25,5.00){\line(1,0){5.00}}
%End Rectangle

\linethickness{0.15mm}
%Rectangle(111.25,0.00)(116.25,5.00)  
\put(111.25,0.00){\line(1,0){5.00}}
\put(111.25,0.00){\line(0,1){5.00}}
\put(116.25,0.00){\line(0,1){5.00}}
\put(111.25,5.00){\line(1,0){5.00}}
%End Rectangle

\linethickness{0.15mm}
%Rectangle(116.25,0.00)(121.25,5.00)  
\put(116.25,0.00){\line(1,0){5.00}}
\put(116.25,0.00){\line(0,1){5.00}}
\put(121.25,0.00){\line(0,1){5.00}}
\put(116.25,5.00){\line(1,0){5.00}}
%End Rectangle

\linethickness{0.15mm}
%Rectangle(121.25,0.00)(126.25,5.00)  
\put(121.25,0.00){\line(1,0){5.00}}
\put(121.25,0.00){\line(0,1){5.00}}
\put(126.25,0.00){\line(0,1){5.00}}
\put(121.25,5.00){\line(1,0){5.00}}
%End Rectangle

\linethickness{0.15mm}
%Rectangle(91.25,5.00)(96.25,10.00)  
\put(91.25,5.00){\line(1,0){5.00}}
\put(91.25,5.00){\line(0,1){5.00}}
\put(96.25,5.00){\line(0,1){5.00}}
\put(91.25,10.00){\line(1,0){5.00}}
%End Rectangle

\linethickness{0.15mm}
%Rectangle(96.25,5.00)(101.25,10.00)  
\put(96.25,5.00){\line(1,0){5.00}}
\put(96.25,5.00){\line(0,1){5.00}}
\put(101.25,5.00){\line(0,1){5.00}}
\put(96.25,10.00){\line(1,0){5.00}}
%End Rectangle

\linethickness{0.15mm}
%Rectangle(101.25,5.00)(106.25,10.00)  
\put(101.25,5.00){\line(1,0){5.00}}
\put(101.25,5.00){\line(0,1){5.00}}
\put(106.25,5.00){\line(0,1){5.00}}
\put(101.25,10.00){\line(1,0){5.00}}
%End Rectangle

\linethickness{0.15mm}
%Rectangle(106.25,5.00)(111.25,10.00)  
\put(106.25,5.00){\line(1,0){5.00}}
\put(106.25,5.00){\line(0,1){5.00}}
\put(111.25,5.00){\line(0,1){5.00}}
\put(106.25,10.00){\line(1,0){5.00}}
%End Rectangle

\linethickness{0.15mm}
%Rectangle(111.25,5.00)(116.25,10.00)  
\put(111.25,5.00){\line(1,0){5.00}}
\put(111.25,5.00){\line(0,1){5.00}}
\put(116.25,5.00){\line(0,1){5.00}}
\put(111.25,10.00){\line(1,0){5.00}}
%End Rectangle

\linethickness{0.15mm}
%Rectangle(116.25,5.00)(121.25,10.00)  
\put(116.25,5.00){\line(1,0){5.00}}
\put(116.25,5.00){\line(0,1){5.00}}
\put(121.25,5.00){\line(0,1){5.00}}
\put(116.25,10.00){\line(1,0){5.00}}
%End Rectangle

\linethickness{0.15mm}
%Rectangle(121.25,5.00)(126.25,10.00)  
\put(121.25,5.00){\line(1,0){5.00}}
\put(121.25,5.00){\line(0,1){5.00}}
\put(126.25,5.00){\line(0,1){5.00}}
\put(121.25,10.00){\line(1,0){5.00}}
%End Rectangle

\linethickness{0.15mm}
%Rectangle(96.25,10.00)(101.25,15.00)  
\put(96.25,10.00){\line(1,0){5.00}}
\put(96.25,10.00){\line(0,1){5.00}}
\put(101.25,10.00){\line(0,1){5.00}}
\put(96.25,15.00){\line(1,0){5.00}}
%End Rectangle

\linethickness{0.15mm}
%Rectangle(101.25,10.00)(106.25,15.00)  
\put(101.25,10.00){\line(1,0){5.00}}
\put(101.25,10.00){\line(0,1){5.00}}
\put(106.25,10.00){\line(0,1){5.00}}
\put(101.25,15.00){\line(1,0){5.00}}
%End Rectangle

\linethickness{0.15mm}
%Rectangle(106.25,10.00)(111.25,15.00)  
\put(106.25,10.00){\line(1,0){5.00}}
\put(106.25,10.00){\line(0,1){5.00}}
\put(111.25,10.00){\line(0,1){5.00}}
\put(106.25,15.00){\line(1,0){5.00}}
%End Rectangle

\linethickness{0.15mm}
%Rectangle(111.25,10.00)(116.25,15.00)  
\put(111.25,10.00){\line(1,0){5.00}}
\put(111.25,10.00){\line(0,1){5.00}}
\put(116.25,10.00){\line(0,1){5.00}}
\put(111.25,15.00){\line(1,0){5.00}}
%End Rectangle

\linethickness{0.15mm}
%Rectangle(116.25,10.00)(121.25,15.00)  
\put(116.25,10.00){\line(1,0){5.00}}
\put(116.25,10.00){\line(0,1){5.00}}
\put(121.25,10.00){\line(0,1){5.00}}
\put(116.25,15.00){\line(1,0){5.00}}
%End Rectangle

\linethickness{0.15mm}
%Rectangle(121.25,10.00)(126.25,15.00)  
\put(121.25,10.00){\line(1,0){5.00}}
\put(121.25,10.00){\line(0,1){5.00}}
\put(126.25,10.00){\line(0,1){5.00}}
\put(121.25,15.00){\line(1,0){5.00}}
%End Rectangle

\linethickness{0.15mm}
%Rectangle(101.25,15.00)(106.25,20.00)  
\put(101.25,15.00){\line(1,0){5.00}}
\put(101.25,15.00){\line(0,1){5.00}}
\put(106.25,15.00){\line(0,1){5.00}}
\put(101.25,20.00){\line(1,0){5.00}}
%End Rectangle

\linethickness{0.15mm}
%Rectangle(106.25,15.00)(111.25,20.00)  
\put(106.25,15.00){\line(1,0){5.00}}
\put(106.25,15.00){\line(0,1){5.00}}
\put(111.25,15.00){\line(0,1){5.00}}
\put(106.25,20.00){\line(1,0){5.00}}
%End Rectangle

\linethickness{0.15mm}
%Rectangle(111.25,15.00)(116.25,20.00)  
\put(111.25,15.00){\line(1,0){5.00}}
\put(111.25,15.00){\line(0,1){5.00}}
\put(116.25,15.00){\line(0,1){5.00}}
\put(111.25,20.00){\line(1,0){5.00}}
%End Rectangle

\linethickness{0.15mm}
%Rectangle(116.25,15.00)(121.25,20.00)  
\put(116.25,15.00){\line(1,0){5.00}}
\put(116.25,15.00){\line(0,1){5.00}}
\put(121.25,15.00){\line(0,1){5.00}}
\put(116.25,20.00){\line(1,0){5.00}}
%End Rectangle

\linethickness{0.15mm}
%Rectangle(106.25,20.00)(111.25,25.00)  
\put(106.25,20.00){\line(1,0){5.00}}
\put(106.25,20.00){\line(0,1){5.00}}
\put(111.25,20.00){\line(0,1){5.00}}
\put(106.25,25.00){\line(1,0){5.00}}
%End Rectangle

\linethickness{0.15mm}
%Rectangle(111.25,20.00)(116.25,25.00)  
\put(111.25,20.00){\line(1,0){5.00}}
\put(111.25,20.00){\line(0,1){5.00}}
\put(116.25,20.00){\line(0,1){5.00}}
\put(111.25,25.00){\line(1,0){5.00}}
%End Rectangle

\linethickness{0.15mm}
%Rectangle(116.25,20.00)(121.25,25.00)  
\put(116.25,20.00){\line(1,0){5.00}}
\put(116.25,20.00){\line(0,1){5.00}}
\put(121.25,20.00){\line(0,1){5.00}}
\put(116.25,25.00){\line(1,0){5.00}}
%End Rectangle

\linethickness{0.15mm}
%Rectangle(111.25,25.00)(116.25,30.00)  
\put(111.25,25.00){\line(1,0){5.00}}
\put(111.25,25.00){\line(0,1){5.00}}
\put(116.25,25.00){\line(0,1){5.00}}
\put(111.25,30.00){\line(1,0){5.00}}
%End Rectangle

\linethickness{0.15mm}
%Rectangle(86.25,0.00)(91.25,5.00)  blacken
\put(86.25,0.00){\rule{5.00\unitlength}{5.00\unitlength}}
%End Rectangle

\linethickness{0.15mm}
%Rectangle(96.25,0.00)(101.25,15.00)  blacken
\put(96.25,0.00){\rule{5.00\unitlength}{15.00\unitlength}}
%End Rectangle

\linethickness{0.15mm}
%Rectangle(106.25,0.00)(111.25,25.00)  blacken
\put(106.25,0.00){\rule{5.00\unitlength}{25.00\unitlength}}
%End Rectangle

\linethickness{0.15mm}
%Rectangle(116.25,0.00)(121.25,25.00)  blacken
\put(116.25,0.00){\rule{5.00\unitlength}{25.00\unitlength}}
%End Rectangle

\linethickness{0.15mm}
%Bezier 0 0(30.00,30.00)(45.00,30.00)(45.00,30.00)
\qbezier(30.00,30.00)(45.00,30.00)(45.00,30.00)
%End Bezier

\linethickness{0.15mm}
%Polygon 0 0(81.25,30.00)(96.25,30.00) 
\put(81.25,30.00){\line(1,0){15.00}}
%End Polygon

\linethickness{0.15mm}
%Polygon 0 0(93.75,31.25)(96.25,30.00) 
\multiput(93.75,31.25)(0.25,-0.13){10}{\line(1,0){0.25}}
%End Polygon

\linethickness{0.15mm}
%Polygon 0 0(96.25,30.00)(93.75,28.75) 
\multiput(93.75,28.75)(0.25,0.13){10}{\line(1,0){0.25}}
%End Polygon

\put(37.50,33.75){\makebox(0,0)[cc]{$(1)$}}

\put(88.75,33.75){\makebox(0,0)[cc]{$(2)$}}

\linethickness{0.15mm}
%Polygon 0 0(43.13,31.25)(45.00,30.00) 
\multiput(43.13,31.25)(0.19,-0.13){10}{\line(1,0){0.19}}
%End Polygon

\linethickness{0.15mm}
%Polygon 0 0(45.00,30.00)(43.13,28.75) 
\multiput(43.13,28.75)(0.19,0.13){10}{\line(1,0){0.19}}
%End Polygon

\end{picture} \\

It is easy to see that these obtained graphs are characterized by the the property that they consist of a finite number of unit squares and regarded from left to right they increase one square at a time untill at some point they are only allowed to be non increasing. The underlying uncoloured graphs are usually called Castelnuovo diagrams or graphs \cite{Davis}. \\

Next we consider the following action on the coloured Castelnuovo graph:
\begin{enumerate}
\item[(3)]
delete one white and black unit square, both on top and on the at most right position as possible.
\end{enumerate}
We repeat (3) as many times as possible in such a way that after every removement the underlying uncoloured graph is a valid Castelnuovo graph. It is easy to see that the inequality \eqref{ineq} holds if it holds after applying (3). We then show that applying (3) a finite number of times we obtain a ``maximal" diagram of the form \\

%\input maximal.tex \\
%Created by jPicEdt 1.x
%Standard LaTeX format (emulated lines)
%Wed Jul 06 22:49:17 CEST 2005
\unitlength 1mm
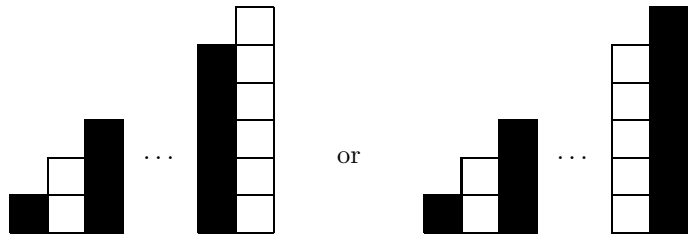
\begin{picture}(110.00,30.00)(0,0)

\linethickness{0.15mm}
%Rectangle(20.00,0.00)(25.00,5.00)  
\put(20.00,0.00){\line(1,0){5.00}}
\put(20.00,0.00){\line(0,1){5.00}}
\put(25.00,0.00){\line(0,1){5.00}}
\put(20.00,5.00){\line(1,0){5.00}}
%End Rectangle

\linethickness{0.15mm}
%Rectangle(25.00,5.00)(30.00,10.00)  
\put(25.00,5.00){\line(1,0){5.00}}
\put(25.00,5.00){\line(0,1){5.00}}
\put(30.00,5.00){\line(0,1){5.00}}
\put(25.00,10.00){\line(1,0){5.00}}
%End Rectangle

\linethickness{0.15mm}
%Rectangle(25.00,0.00)(30.00,5.00)  
\put(25.00,0.00){\line(1,0){5.00}}
\put(25.00,0.00){\line(0,1){5.00}}
\put(30.00,0.00){\line(0,1){5.00}}
\put(25.00,5.00){\line(1,0){5.00}}
%End Rectangle

\linethickness{0.15mm}
%Rectangle(30.00,10.00)(35.00,15.00)  
\put(30.00,10.00){\line(1,0){5.00}}
\put(30.00,10.00){\line(0,1){5.00}}
\put(35.00,10.00){\line(0,1){5.00}}
\put(30.00,15.00){\line(1,0){5.00}}
%End Rectangle

\linethickness{0.15mm}
%Rectangle(30.00,5.00)(35.00,10.00)  
\put(30.00,5.00){\line(1,0){5.00}}
\put(30.00,5.00){\line(0,1){5.00}}
\put(35.00,5.00){\line(0,1){5.00}}
\put(30.00,10.00){\line(1,0){5.00}}
%End Rectangle

\linethickness{0.15mm}
%Rectangle(30.00,0.00)(35.00,5.00)  
\put(30.00,0.00){\line(1,0){5.00}}
\put(30.00,0.00){\line(0,1){5.00}}
\put(35.00,0.00){\line(0,1){5.00}}
\put(30.00,5.00){\line(1,0){5.00}}
%End Rectangle

\linethickness{0.15mm}
%Rectangle(45.00,20.00)(50.00,25.00)  
\put(45.00,20.00){\line(1,0){5.00}}
\put(45.00,20.00){\line(0,1){5.00}}
\put(50.00,20.00){\line(0,1){5.00}}
\put(45.00,25.00){\line(1,0){5.00}}
%End Rectangle

\linethickness{0.15mm}
%Rectangle(45.00,15.00)(50.00,20.00)  
\put(45.00,15.00){\line(1,0){5.00}}
\put(45.00,15.00){\line(0,1){5.00}}
\put(50.00,15.00){\line(0,1){5.00}}
\put(45.00,20.00){\line(1,0){5.00}}
%End Rectangle

\linethickness{0.15mm}
%Rectangle(45.00,10.00)(50.00,15.00)  
\put(45.00,10.00){\line(1,0){5.00}}
\put(45.00,10.00){\line(0,1){5.00}}
\put(50.00,10.00){\line(0,1){5.00}}
\put(45.00,15.00){\line(1,0){5.00}}
%End Rectangle

\linethickness{0.15mm}
%Rectangle(45.00,5.00)(50.00,10.00)  
\put(45.00,5.00){\line(1,0){5.00}}
\put(45.00,5.00){\line(0,1){5.00}}
\put(50.00,5.00){\line(0,1){5.00}}
\put(45.00,10.00){\line(1,0){5.00}}
%End Rectangle

\linethickness{0.15mm}
%Rectangle(45.00,0.00)(50.00,5.00)  
\put(45.00,0.00){\line(1,0){5.00}}
\put(45.00,0.00){\line(0,1){5.00}}
\put(50.00,0.00){\line(0,1){5.00}}
\put(45.00,5.00){\line(1,0){5.00}}
%End Rectangle

\linethickness{0.15mm}
%Rectangle(50.00,25.00)(55.00,30.00)  
\put(50.00,25.00){\line(1,0){5.00}}
\put(50.00,25.00){\line(0,1){5.00}}
\put(55.00,25.00){\line(0,1){5.00}}
\put(50.00,30.00){\line(1,0){5.00}}
%End Rectangle

\linethickness{0.15mm}
%Rectangle(50.00,20.00)(55.00,25.00)  
\put(50.00,20.00){\line(1,0){5.00}}
\put(50.00,20.00){\line(0,1){5.00}}
\put(55.00,20.00){\line(0,1){5.00}}
\put(50.00,25.00){\line(1,0){5.00}}
%End Rectangle

\linethickness{0.15mm}
%Rectangle(50.00,15.00)(55.00,20.00)  
\put(50.00,15.00){\line(1,0){5.00}}
\put(50.00,15.00){\line(0,1){5.00}}
\put(55.00,15.00){\line(0,1){5.00}}
\put(50.00,20.00){\line(1,0){5.00}}
%End Rectangle

\linethickness{0.15mm}
%Rectangle(50.00,10.00)(55.00,15.00)  
\put(50.00,10.00){\line(1,0){5.00}}
\put(50.00,10.00){\line(0,1){5.00}}
\put(55.00,10.00){\line(0,1){5.00}}
\put(50.00,15.00){\line(1,0){5.00}}
%End Rectangle

\linethickness{0.15mm}
%Rectangle(50.00,5.00)(55.00,10.00)  
\put(50.00,5.00){\line(1,0){5.00}}
\put(50.00,5.00){\line(0,1){5.00}}
\put(55.00,5.00){\line(0,1){5.00}}
\put(50.00,10.00){\line(1,0){5.00}}
%End Rectangle

\linethickness{0.15mm}
%Rectangle(50.00,0.00)(55.00,5.00)  
\put(50.00,0.00){\line(1,0){5.00}}
\put(50.00,0.00){\line(0,1){5.00}}
\put(55.00,0.00){\line(0,1){5.00}}
\put(50.00,5.00){\line(1,0){5.00}}
%End Rectangle

\put(40.00,10.00){\makebox(0,0)[cc]{$\dots$}}

\linethickness{0.15mm}
%Rectangle(75.00,0.00)(80.00,5.00)  
\put(75.00,0.00){\line(1,0){5.00}}
\put(75.00,0.00){\line(0,1){5.00}}
\put(80.00,0.00){\line(0,1){5.00}}
\put(75.00,5.00){\line(1,0){5.00}}
%End Rectangle

\linethickness{0.15mm}
%Rectangle(80.00,5.00)(85.00,10.00)  
\put(80.00,5.00){\line(1,0){5.00}}
\put(80.00,5.00){\line(0,1){5.00}}
\put(85.00,5.00){\line(0,1){5.00}}
\put(80.00,10.00){\line(1,0){5.00}}
%End Rectangle

\linethickness{0.15mm}
%Rectangle(80.00,0.00)(85.00,5.00)  
\put(80.00,0.00){\line(1,0){5.00}}
\put(80.00,0.00){\line(0,1){5.00}}
\put(85.00,0.00){\line(0,1){5.00}}
\put(80.00,5.00){\line(1,0){5.00}}
%End Rectangle

\linethickness{0.15mm}
%Rectangle(85.00,10.00)(90.00,15.00)  
\put(85.00,10.00){\line(1,0){5.00}}
\put(85.00,10.00){\line(0,1){5.00}}
\put(90.00,10.00){\line(0,1){5.00}}
\put(85.00,15.00){\line(1,0){5.00}}
%End Rectangle

\linethickness{0.15mm}
%Rectangle(85.00,5.00)(90.00,10.00)  
\put(85.00,5.00){\line(1,0){5.00}}
\put(85.00,5.00){\line(0,1){5.00}}
\put(90.00,5.00){\line(0,1){5.00}}
\put(85.00,10.00){\line(1,0){5.00}}
%End Rectangle

\linethickness{0.15mm}
%Rectangle(85.00,0.00)(90.00,5.00)  
\put(85.00,0.00){\line(1,0){5.00}}
\put(85.00,0.00){\line(0,1){5.00}}
\put(90.00,0.00){\line(0,1){5.00}}
\put(85.00,5.00){\line(1,0){5.00}}
%End Rectangle

\linethickness{0.15mm}
%Rectangle(100.00,20.00)(105.00,25.00)  
\put(100.00,20.00){\line(1,0){5.00}}
\put(100.00,20.00){\line(0,1){5.00}}
\put(105.00,20.00){\line(0,1){5.00}}
\put(100.00,25.00){\line(1,0){5.00}}
%End Rectangle

\linethickness{0.15mm}
%Rectangle(100.00,15.00)(105.00,20.00)  
\put(100.00,15.00){\line(1,0){5.00}}
\put(100.00,15.00){\line(0,1){5.00}}
\put(105.00,15.00){\line(0,1){5.00}}
\put(100.00,20.00){\line(1,0){5.00}}
%End Rectangle

\linethickness{0.15mm}
%Rectangle(100.00,10.00)(105.00,15.00)  
\put(100.00,10.00){\line(1,0){5.00}}
\put(100.00,10.00){\line(0,1){5.00}}
\put(105.00,10.00){\line(0,1){5.00}}
\put(100.00,15.00){\line(1,0){5.00}}
%End Rectangle

\linethickness{0.15mm}
%Rectangle(100.00,5.00)(105.00,10.00)  
\put(100.00,5.00){\line(1,0){5.00}}
\put(100.00,5.00){\line(0,1){5.00}}
\put(105.00,5.00){\line(0,1){5.00}}
\put(100.00,10.00){\line(1,0){5.00}}
%End Rectangle

\linethickness{0.15mm}
%Rectangle(100.00,0.00)(105.00,5.00)  
\put(100.00,0.00){\line(1,0){5.00}}
\put(100.00,0.00){\line(0,1){5.00}}
\put(105.00,0.00){\line(0,1){5.00}}
\put(100.00,5.00){\line(1,0){5.00}}
%End Rectangle

\linethickness{0.15mm}
%Rectangle(105.00,25.00)(110.00,30.00)  
\put(105.00,25.00){\line(1,0){5.00}}
\put(105.00,25.00){\line(0,1){5.00}}
\put(110.00,25.00){\line(0,1){5.00}}
\put(105.00,30.00){\line(1,0){5.00}}
%End Rectangle

\linethickness{0.15mm}
%Rectangle(105.00,20.00)(110.00,25.00)  
\put(105.00,20.00){\line(1,0){5.00}}
\put(105.00,20.00){\line(0,1){5.00}}
\put(110.00,20.00){\line(0,1){5.00}}
\put(105.00,25.00){\line(1,0){5.00}}
%End Rectangle

\linethickness{0.15mm}
%Rectangle(105.00,15.00)(110.00,20.00)  
\put(105.00,15.00){\line(1,0){5.00}}
\put(105.00,15.00){\line(0,1){5.00}}
\put(110.00,15.00){\line(0,1){5.00}}
\put(105.00,20.00){\line(1,0){5.00}}
%End Rectangle

\linethickness{0.15mm}
%Rectangle(105.00,10.00)(110.00,15.00)  
\put(105.00,10.00){\line(1,0){5.00}}
\put(105.00,10.00){\line(0,1){5.00}}
\put(110.00,10.00){\line(0,1){5.00}}
\put(105.00,15.00){\line(1,0){5.00}}
%End Rectangle

\linethickness{0.15mm}
%Rectangle(105.00,5.00)(110.00,10.00)  
\put(105.00,5.00){\line(1,0){5.00}}
\put(105.00,5.00){\line(0,1){5.00}}
\put(110.00,5.00){\line(0,1){5.00}}
\put(105.00,10.00){\line(1,0){5.00}}
%End Rectangle

\linethickness{0.15mm}
%Rectangle(105.00,0.00)(110.00,5.00)  
\put(105.00,0.00){\line(1,0){5.00}}
\put(105.00,0.00){\line(0,1){5.00}}
\put(110.00,0.00){\line(0,1){5.00}}
\put(105.00,5.00){\line(1,0){5.00}}
%End Rectangle

\put(95.00,10.00){\makebox(0,0)[cc]{$\dots$}}

\linethickness{0.15mm}
%Rectangle(20.00,0.00)(25.00,5.00)  blacken
\put(20.00,0.00){\rule{5.00\unitlength}{5.00\unitlength}}
%End Rectangle

\linethickness{0.15mm}
%Rectangle(30.00,0.00)(35.00,15.00)  blacken
\put(30.00,0.00){\rule{5.00\unitlength}{15.00\unitlength}}
%End Rectangle

\linethickness{0.15mm}
%Rectangle(45.00,0.00)(50.00,25.00)  blacken
\put(45.00,0.00){\rule{5.00\unitlength}{25.00\unitlength}}
%End Rectangle

\linethickness{0.15mm}
%Rectangle(75.00,0.00)(80.00,5.00)  blacken
\put(75.00,0.00){\rule{5.00\unitlength}{5.00\unitlength}}
%End Rectangle

\linethickness{0.15mm}
%Rectangle(85.00,0.00)(90.00,15.00)  blacken
\put(85.00,0.00){\rule{5.00\unitlength}{15.00\unitlength}}
%End Rectangle

\linethickness{0.15mm}
%Rectangle(105.00,0.00)(110.00,30.00)  blacken
\put(105.00,0.00){\rule{5.00\unitlength}{30.00\unitlength}}
%End Rectangle

\put(65.00,10.00){\makebox(0,0)[cc]{or}}

\end{picture} \\

for which \eqref{ineq} is (trivially) true. This proves that the condition \eqref{ineq} is necessary. To prove that \eqref{ineq} is sufficient we show that there exists a (coloured) Castelnuovo graph of the form \\

%\input maximalcastel.tex 
%Created by jPicEdt 1.x
%Standard LaTeX format (emulated lines)
%Thu Jul 07 17:59:09 CEST 2005
\unitlength 1mm
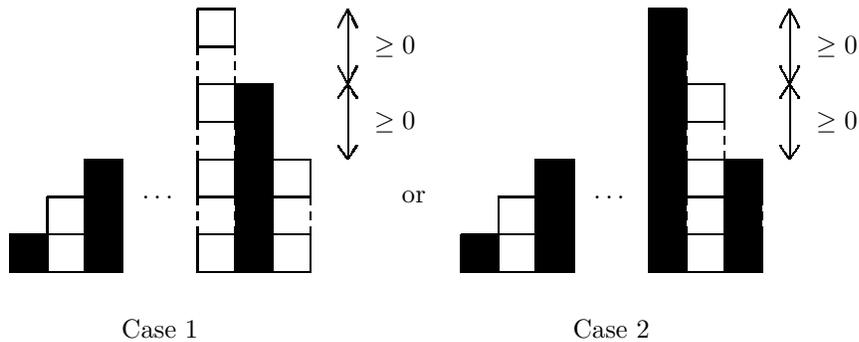
\begin{picture}(115.00,45.00)(0,0)

\linethickness{0.15mm}
%Rectangle(65.00,10.00)(70.00,15.00)  
\put(65.00,10.00){\line(1,0){5.00}}
\put(65.00,10.00){\line(0,1){5.00}}
\put(70.00,10.00){\line(0,1){5.00}}
\put(65.00,15.00){\line(1,0){5.00}}
%End Rectangle

\linethickness{0.15mm}
%Rectangle(70.00,15.00)(75.00,20.00)  
\put(70.00,15.00){\line(1,0){5.00}}
\put(70.00,15.00){\line(0,1){5.00}}
\put(75.00,15.00){\line(0,1){5.00}}
\put(70.00,20.00){\line(1,0){5.00}}
%End Rectangle

\linethickness{0.15mm}
%Rectangle(70.00,10.00)(75.00,15.00)  
\put(70.00,10.00){\line(1,0){5.00}}
\put(70.00,10.00){\line(0,1){5.00}}
\put(75.00,10.00){\line(0,1){5.00}}
\put(70.00,15.00){\line(1,0){5.00}}
%End Rectangle

\linethickness{0.15mm}
%Rectangle(75.00,20.00)(80.00,25.00)  
\put(75.00,20.00){\line(1,0){5.00}}
\put(75.00,20.00){\line(0,1){5.00}}
\put(80.00,20.00){\line(0,1){5.00}}
\put(75.00,25.00){\line(1,0){5.00}}
%End Rectangle

\linethickness{0.15mm}
%Rectangle(75.00,15.00)(80.00,20.00)  
\put(75.00,15.00){\line(1,0){5.00}}
\put(75.00,15.00){\line(0,1){5.00}}
\put(80.00,15.00){\line(0,1){5.00}}
\put(75.00,20.00){\line(1,0){5.00}}
%End Rectangle

\linethickness{0.15mm}
%Rectangle(75.00,10.00)(80.00,15.00)  
\put(75.00,10.00){\line(1,0){5.00}}
\put(75.00,10.00){\line(0,1){5.00}}
\put(80.00,10.00){\line(0,1){5.00}}
\put(75.00,15.00){\line(1,0){5.00}}
%End Rectangle

\linethickness{0.15mm}
%Rectangle(90.00,30.00)(95.00,35.00)  
\put(90.00,30.00){\line(1,0){5.00}}
\put(90.00,30.00){\line(0,1){5.00}}
\put(95.00,30.00){\line(0,1){5.00}}
\put(90.00,35.00){\line(1,0){5.00}}
%End Rectangle

\linethickness{0.15mm}
%Rectangle(90.00,20.00)(95.00,25.00)  
\put(90.00,20.00){\line(1,0){5.00}}
\put(90.00,20.00){\line(0,1){5.00}}
\put(95.00,20.00){\line(0,1){5.00}}
\put(90.00,25.00){\line(1,0){5.00}}
%End Rectangle

\linethickness{0.15mm}
%Rectangle(90.00,10.00)(95.00,15.00)  
\put(90.00,10.00){\line(1,0){5.00}}
\put(90.00,10.00){\line(0,1){5.00}}
\put(95.00,10.00){\line(0,1){5.00}}
\put(90.00,15.00){\line(1,0){5.00}}
%End Rectangle

\linethickness{0.15mm}
%Rectangle(90.00,40.00)(95.00,45.00)  
\put(90.00,40.00){\line(1,0){5.00}}
\put(90.00,40.00){\line(0,1){5.00}}
\put(95.00,40.00){\line(0,1){5.00}}
\put(90.00,45.00){\line(1,0){5.00}}
%End Rectangle

\linethickness{0.15mm}
%Rectangle(100.00,10.00)(105.00,15.00)  
\put(100.00,10.00){\line(1,0){5.00}}
\put(100.00,10.00){\line(0,1){5.00}}
\put(105.00,10.00){\line(0,1){5.00}}
\put(100.00,15.00){\line(1,0){5.00}}
%End Rectangle

\linethickness{0.15mm}
%Rectangle(100.00,20.00)(105.00,25.00)  
\put(100.00,20.00){\line(1,0){5.00}}
\put(100.00,20.00){\line(0,1){5.00}}
\put(105.00,20.00){\line(0,1){5.00}}
\put(100.00,25.00){\line(1,0){5.00}}
%End Rectangle

\linethickness{0.15mm}
%Rectangle(95.00,30.00)(100.00,35.00)  
\put(95.00,30.00){\line(1,0){5.00}}
\put(95.00,30.00){\line(0,1){5.00}}
\put(100.00,30.00){\line(0,1){5.00}}
\put(95.00,35.00){\line(1,0){5.00}}
%End Rectangle

\linethickness{0.15mm}
%Rectangle(95.00,20.00)(100.00,25.00)  
\put(95.00,20.00){\line(1,0){5.00}}
\put(95.00,20.00){\line(0,1){5.00}}
\put(100.00,20.00){\line(0,1){5.00}}
\put(95.00,25.00){\line(1,0){5.00}}
%End Rectangle

\linethickness{0.15mm}
%Rectangle(95.00,10.00)(100.00,15.00)  
\put(95.00,10.00){\line(1,0){5.00}}
\put(95.00,10.00){\line(0,1){5.00}}
\put(100.00,10.00){\line(0,1){5.00}}
\put(95.00,15.00){\line(1,0){5.00}}
%End Rectangle

\linethickness{0.15mm}
%Polygon 0 0(90.00,40.00)(90.00,15.00) dash=1.00
\multiput(90.00,15.00)(0,2.00){13}{\line(0,1){1.00}}
%End Polygon

\linethickness{0.15mm}
%Polygon 0 0(95.00,40.00)(95.00,15.00) dash=1.00
\multiput(95.00,15.00)(0,2.00){13}{\line(0,1){1.00}}
%End Polygon

\linethickness{0.15mm}
%Polygon 0 0(100.00,30.00)(100.00,15.00) dash=1.00
\multiput(100.00,15.00)(0,2.00){8}{\line(0,1){1.00}}
%End Polygon

\linethickness{0.15mm}
%Polygon 0 0(105.00,20.00)(105.00,15.00) dash=1.00
\multiput(105.00,15.00)(0,2.00){3}{\line(0,1){1.00}}
%End Polygon

\put(85.00,20.00){\makebox(0,0)[cc]{$\dots$}}

\linethickness{0.15mm}
%Rectangle(65.00,10.00)(70.00,15.00)  blacken
\put(65.00,10.00){\rule{5.00\unitlength}{5.00\unitlength}}
%End Rectangle

\linethickness{0.15mm}
%Rectangle(75.00,10.00)(80.00,25.00)  blacken
\put(75.00,10.00){\rule{5.00\unitlength}{15.00\unitlength}}
%End Rectangle

\linethickness{0.15mm}
%Rectangle(5.00,10.00)(10.00,15.00)  
\put(5.00,10.00){\line(1,0){5.00}}
\put(5.00,10.00){\line(0,1){5.00}}
\put(10.00,10.00){\line(0,1){5.00}}
\put(5.00,15.00){\line(1,0){5.00}}
%End Rectangle

\linethickness{0.15mm}
%Rectangle(10.00,15.00)(15.00,20.00)  
\put(10.00,15.00){\line(1,0){5.00}}
\put(10.00,15.00){\line(0,1){5.00}}
\put(15.00,15.00){\line(0,1){5.00}}
\put(10.00,20.00){\line(1,0){5.00}}
%End Rectangle

\linethickness{0.15mm}
%Rectangle(10.00,10.00)(15.00,15.00)  
\put(10.00,10.00){\line(1,0){5.00}}
\put(10.00,10.00){\line(0,1){5.00}}
\put(15.00,10.00){\line(0,1){5.00}}
\put(10.00,15.00){\line(1,0){5.00}}
%End Rectangle

\linethickness{0.15mm}
%Rectangle(15.00,20.00)(20.00,25.00)  
\put(15.00,20.00){\line(1,0){5.00}}
\put(15.00,20.00){\line(0,1){5.00}}
\put(20.00,20.00){\line(0,1){5.00}}
\put(15.00,25.00){\line(1,0){5.00}}
%End Rectangle

\linethickness{0.15mm}
%Rectangle(15.00,15.00)(20.00,20.00)  
\put(15.00,15.00){\line(1,0){5.00}}
\put(15.00,15.00){\line(0,1){5.00}}
\put(20.00,15.00){\line(0,1){5.00}}
\put(15.00,20.00){\line(1,0){5.00}}
%End Rectangle

\linethickness{0.15mm}
%Rectangle(15.00,10.00)(20.00,15.00)  
\put(15.00,10.00){\line(1,0){5.00}}
\put(15.00,10.00){\line(0,1){5.00}}
\put(20.00,10.00){\line(0,1){5.00}}
\put(15.00,15.00){\line(1,0){5.00}}
%End Rectangle

\linethickness{0.15mm}
%Rectangle(30.00,30.00)(35.00,35.00)  
\put(30.00,30.00){\line(1,0){5.00}}
\put(30.00,30.00){\line(0,1){5.00}}
\put(35.00,30.00){\line(0,1){5.00}}
\put(30.00,35.00){\line(1,0){5.00}}
%End Rectangle

\linethickness{0.15mm}
%Rectangle(30.00,20.00)(35.00,25.00)  
\put(30.00,20.00){\line(1,0){5.00}}
\put(30.00,20.00){\line(0,1){5.00}}
\put(35.00,20.00){\line(0,1){5.00}}
\put(30.00,25.00){\line(1,0){5.00}}
%End Rectangle

\linethickness{0.15mm}
%Rectangle(30.00,10.00)(35.00,15.00)  
\put(30.00,10.00){\line(1,0){5.00}}
\put(30.00,10.00){\line(0,1){5.00}}
\put(35.00,10.00){\line(0,1){5.00}}
\put(30.00,15.00){\line(1,0){5.00}}
%End Rectangle

\linethickness{0.15mm}
%Rectangle(30.00,40.00)(35.00,45.00)  
\put(30.00,40.00){\line(1,0){5.00}}
\put(30.00,40.00){\line(0,1){5.00}}
\put(35.00,40.00){\line(0,1){5.00}}
\put(30.00,45.00){\line(1,0){5.00}}
%End Rectangle

\linethickness{0.15mm}
%Rectangle(40.00,10.00)(45.00,15.00)  
\put(40.00,10.00){\line(1,0){5.00}}
\put(40.00,10.00){\line(0,1){5.00}}
\put(45.00,10.00){\line(0,1){5.00}}
\put(40.00,15.00){\line(1,0){5.00}}
%End Rectangle

\linethickness{0.15mm}
%Rectangle(40.00,20.00)(45.00,25.00)  
\put(40.00,20.00){\line(1,0){5.00}}
\put(40.00,20.00){\line(0,1){5.00}}
\put(45.00,20.00){\line(0,1){5.00}}
\put(40.00,25.00){\line(1,0){5.00}}
%End Rectangle

\linethickness{0.15mm}
%Rectangle(35.00,30.00)(40.00,35.00)  
\put(35.00,30.00){\line(1,0){5.00}}
\put(35.00,30.00){\line(0,1){5.00}}
\put(40.00,30.00){\line(0,1){5.00}}
\put(35.00,35.00){\line(1,0){5.00}}
%End Rectangle

\linethickness{0.15mm}
%Rectangle(35.00,20.00)(40.00,25.00)  
\put(35.00,20.00){\line(1,0){5.00}}
\put(35.00,20.00){\line(0,1){5.00}}
\put(40.00,20.00){\line(0,1){5.00}}
\put(35.00,25.00){\line(1,0){5.00}}
%End Rectangle

\linethickness{0.15mm}
%Rectangle(35.00,10.00)(40.00,15.00)  
\put(35.00,10.00){\line(1,0){5.00}}
\put(35.00,10.00){\line(0,1){5.00}}
\put(40.00,10.00){\line(0,1){5.00}}
\put(35.00,15.00){\line(1,0){5.00}}
%End Rectangle

\linethickness{0.15mm}
%Polygon 0 0(30.00,40.00)(30.00,15.00) dash=1.00
\multiput(30.00,15.00)(0,2.00){13}{\line(0,1){1.00}}
%End Polygon

\linethickness{0.15mm}
%Polygon 0 0(35.00,40.00)(35.00,15.00) dash=1.00
\multiput(35.00,15.00)(0,2.00){13}{\line(0,1){1.00}}
%End Polygon

\linethickness{0.15mm}
%Polygon 0 0(40.00,30.00)(40.00,15.00) dash=1.00
\multiput(40.00,15.00)(0,2.00){8}{\line(0,1){1.00}}
%End Polygon

\linethickness{0.15mm}
%Polygon 0 0(45.00,20.00)(45.00,15.00) dash=1.00
\multiput(45.00,15.00)(0,2.00){3}{\line(0,1){1.00}}
%End Polygon

\put(25.00,20.00){\makebox(0,0)[cc]{$\dots$}}

\linethickness{0.15mm}
%Rectangle(5.00,10.00)(10.00,15.00)  blacken
\put(5.00,10.00){\rule{5.00\unitlength}{5.00\unitlength}}
%End Rectangle

\linethickness{0.15mm}
%Rectangle(15.00,10.00)(20.00,25.00)  blacken
\put(15.00,10.00){\rule{5.00\unitlength}{15.00\unitlength}}
%End Rectangle

\put(58.75,20.00){\makebox(0,0)[cc]{or}}

\linethickness{0.15mm}
%Rectangle(90.00,10.00)(95.00,45.00)  blacken
\put(90.00,10.00){\rule{5.00\unitlength}{35.00\unitlength}}
%End Rectangle

\linethickness{0.15mm}
%Rectangle(100.00,10.00)(105.00,25.00)  blacken
\put(100.00,10.00){\rule{5.00\unitlength}{15.00\unitlength}}
%End Rectangle

\linethickness{0.15mm}
%Rectangle(35.00,10.00)(40.00,35.00)  blacken
\put(35.00,10.00){\rule{5.00\unitlength}{25.00\unitlength}}
%End Rectangle

\linethickness{0.15mm}
%Polygon 0 0(108.75,45.00)(108.75,35.00) blacken
\put(108.75,35.00){\line(0,1){10.00}}
%End Polygon

\linethickness{0.15mm}
%Polygon 0 0(108.75,35.00)(108.75,25.00) blacken
\put(108.75,25.00){\line(0,1){10.00}}
%End Polygon

\put(115.00,40.00){\makebox(0,0)[cc]{$\geq 0$}}

\put(115.00,30.00){\makebox(0,0)[cc]{$\geq 0$}}

\linethickness{0.15mm}
%Polygon 0 0(108.75,45.00)(107.50,42.50) blacken
\multiput(107.50,42.50)(0.13,0.25){10}{\line(0,1){0.25}}
%End Polygon

\linethickness{0.15mm}
%Polygon 0 0(108.75,45.00)(110.00,42.50) blacken
\multiput(108.75,45.00)(0.13,-0.25){10}{\line(0,-1){0.25}}
%End Polygon

\linethickness{0.15mm}
%Polygon 0 0(108.75,35.00)(107.50,37.50) blacken
\multiput(107.50,37.50)(0.13,-0.25){10}{\line(0,-1){0.25}}
%End Polygon

\linethickness{0.15mm}
%Polygon 0 0(108.75,35.00)(110.00,37.50) blacken
\multiput(108.75,35.00)(0.13,0.25){10}{\line(0,1){0.25}}
%End Polygon

\linethickness{0.15mm}
%Polygon 0 0(108.75,35.00)(107.50,33.13) blacken
\multiput(107.50,33.13)(0.13,0.19){10}{\line(0,1){0.19}}
%End Polygon

\linethickness{0.15mm}
%Polygon 0 0(108.75,35.00)(108.75,35.00) blacken

%End Polygon

\linethickness{0.15mm}
%Polygon 0 0(108.75,35.00)(110.00,33.13) blacken
\multiput(108.75,35.00)(0.13,-0.19){10}{\line(0,-1){0.19}}
%End Polygon

\linethickness{0.15mm}
%Polygon 0 0(107.50,26.88)(108.75,25.00) blacken
\multiput(107.50,26.88)(0.13,-0.19){10}{\line(0,-1){0.19}}
%End Polygon

\linethickness{0.15mm}
%Polygon 0 0(108.75,25.00)(110.00,26.88) blacken
\multiput(108.75,25.00)(0.13,0.19){10}{\line(0,1){0.19}}
%End Polygon

\linethickness{0.15mm}
%Polygon 0 0(50.00,45.00)(50.00,35.00) blacken
\put(50.00,35.00){\line(0,1){10.00}}
%End Polygon

\linethickness{0.15mm}
%Polygon 0 0(50.00,35.00)(50.00,25.00) blacken
\put(50.00,25.00){\line(0,1){10.00}}
%End Polygon

\put(56.25,40.00){\makebox(0,0)[cc]{$\geq 0$}}

\put(56.25,30.00){\makebox(0,0)[cc]{$\geq 0$}}

\linethickness{0.15mm}
%Polygon 0 0(50.00,45.00)(48.75,42.50) blacken
\multiput(48.75,42.50)(0.13,0.25){10}{\line(0,1){0.25}}
%End Polygon

\linethickness{0.15mm}
%Polygon 0 0(50.00,45.00)(51.25,42.50) blacken
\multiput(50.00,45.00)(0.13,-0.25){10}{\line(0,-1){0.25}}
%End Polygon

\linethickness{0.15mm}
%Polygon 0 0(50.00,35.00)(48.75,37.50) blacken
\multiput(48.75,37.50)(0.13,-0.25){10}{\line(0,-1){0.25}}
%End Polygon

\linethickness{0.15mm}
%Polygon 0 0(50.00,35.00)(51.25,37.50) blacken
\multiput(50.00,35.00)(0.13,0.25){10}{\line(0,1){0.25}}
%End Polygon

\linethickness{0.15mm}
%Polygon 0 0(50.00,35.00)(48.75,33.12) blacken
\multiput(48.75,33.12)(0.13,0.19){10}{\line(0,1){0.19}}
%End Polygon

\linethickness{0.15mm}
%Polygon 0 0(50.00,35.00)(50.00,35.00) blacken

%End Polygon

\linethickness{0.15mm}
%Polygon 0 0(50.00,35.00)(51.25,33.12) blacken
\multiput(50.00,35.00)(0.13,-0.19){10}{\line(0,-1){0.19}}
%End Polygon

\linethickness{0.15mm}
%Polygon 0 0(48.75,26.87)(50.00,25.00) blacken
\multiput(48.75,26.87)(0.13,-0.19){10}{\line(0,-1){0.19}}
%End Polygon

\linethickness{0.15mm}
%Polygon 0 0(50.00,25.00)(51.25,26.87) blacken
\multiput(50.00,25.00)(0.13,0.19){10}{\line(0,1){0.19}}
%End Polygon

\put(85.00,2.50){\makebox(0,0)[cc]{Case 2}}

\put(25.00,2.50){\makebox(0,0)[cc]{Case 1}}

\end{picture} 
\newpage
where the sum of black (resp. white) unit squares is equal to $b$ (resp. $w$). By reversing the above proces we find a partition $\la$ for which $(b(\la),w(\la)) = (b,w)$. As a refinement, this partition has distinct parts.
\begin{remark}
The authors found the inequality \eqref{ineq} in Theorem A while investigating Hilbert series of reflexive rank one modules over cubic Artin-Schelter regular $k$-algebras $A$ of global dimension three \cite{AS, ATV1, ATV2}. In this context $k$ is an algebraically closed field of characteristic zero. These graded algebras $A$ are regarded as noncommutative analogues of the coordinate ring of a quadric in $\PP^{3}$. Let us sketch briefly how we obtained \eqref{ineq}. See \cite{DM} for more details. Assume that $A$ is such a cubic algebra. For any reflexive rank one module $M$ over $A$ the Hilbert series of $M$ is (up to shift of grading) of the form 
\begin{equation} \label{Hilb}
h_{M}(t) = h_{A}(t) - \frac{s(t)}{1-t^{2}} + f(t)
\end{equation}
for some $f(t),s(t) \in \ZZ[t,t^{-1}]$. It turns out that $s(t)$ is the generating function of a Castelnuovo function (related to a Castelnuovo diagram, see \S\ref{2.2} for its definition). Writing $(b,w) =(b(s),w(s))$ the equation \eqref{Hilb} implies
\begin{align*}
\dim_{k}A_{l} - \dim_{k}M_{l} = 
\left\{
\begin{array}{ll}
b & \text{ if $l \gg 0$ is even} \\
w & \text{ if $l \gg 0$ is odd}
\end{array}
\right.
\end{align*}
Moreover, if the algebra $A$ is generic then for any Castelnuovo function $s$ there exists a reflexive rank one module $M$ such that (after shift of grading) such that \eqref{Hilb} holds. On the other hand we find $\dim_{k}\Ext^{1}(\Mscr,\Mscr) = 2(b - (b-w)^{2})$ where $\Mscr = \pi M$ is the quotient of $M$ by the maximal finite dimensional %(as a $k$-vector space) 
submodule of $M$. Since this dimension has to be positive we therefore conclude that for any Castelnuovo function $s$ (and hence for any partition $\la$) the inequality \eqref{ineq} holds. 
\end{remark}

The rest of this note is organized as follows. In Section 2 we have included some preliminaries on partitions and Castelnuovo function. We develop their relation which we will need later on. In Section 3 the proof of Theorem A is given. Section 4 presents the proof of Theorem B. Finally in Section 5 we make the connection to \cite[Problem 10]{competition}.

\section{Generalities}

It this section we recall some basic notions. We refer to \cite{Andrews} for an introduction to the theory of partitions. 

\subsection{Partitions and chess Ferrers graphs} \label{2.1}

A {\em partition} $\la$ of a positive integer $n$ is a finite sequence of positive integers $\lambda_{1},\lambda_{2},\dots, \lambda_{r}$ for which
\[
\la_{1} \geq \la_{2} \geq \dots \geq \la_{r} \geq 0  \quad \text{ and } \quad \sum_{i = 1}^{r}\la_{i} = n.
\]
We will often not specify the integer $n$, and put $\la_{i} = 0$ for $i<1$ and $i>r$. %Thus by a partition we simply mean a partition of $n$ for some positive integer $n$. 
The partition $(\la_{1},\la_{2},\dots,\la_{r})$ will be denoted by $\la$ and for convenience we assume that the appearing entries in $\la$ are nonzero. Thus the empty sequence $\la = (\,\,)$ forms the only partition of zero. We refer to the integers $\la_{1}, \dots, \la_{r}$ as the {\em parts} of $\la$. In case all parts of $\la$ are distinct we 
say that $\la$ is a {\em partition in distinct parts}. The sum 
%the appearing integer $r$ as the {\em length} of $\la$ and 
$n = \la_{1} + \la_{2} + \dots + \la_{r}$ is called the {\em weight} of $\la$. Write $\Pscr$ for the set of all partitions (of weight $n$ where $n$ runs through all positive integers). Similary we let $\Dscr \subset \Pscr$ be the set of all partitions in distinct parts. \\

If $\la \in \Pscr$ is a partition we may define a new partition $\la' = (\la'_{1},\la'_{2},\dots,\la'_{r'})$ by defining $\la'_{i}$ as the number of parts of $\la$ that are greater or equal than $i$ (for $i \geq 1$)
\[
\la'_{i} = \# \{ j \mid \la_{j} \geq i \}
\]
The partition $\la'$ is called the {\em conjugate} of $\la$. 
Note that $\weight \la = \weight \la'$.  
%We write $p_{n}(t)$ for its associated polynomial
%\[
%p_{n}(t) = \la_{0} + \la_{1}t + \la_{2}t^{2} + \dots + \la_{l}t^{l-1}
%\]
It is standard to visualize a partition $\la \in \Pscr$ using the graph of the staircase function 
\[ 
F(\la): \RR \r \NN: x \mapsto \la'_{\lfloor x \rfloor} 
\] 
where $\lfloor x \rfloor$ stands for the greatest integer less or equal than $x \in \RR$. We divide the area under this graph $F(\la)$ in unit cases. This graph is called the {\em Ferrers graph} of $\la$. Note that the number of unit squares in the diagram is equal to the weight of $\la$. We label the columns from left to right, and rows from down to up, starting by index number zero.
\begin{example}
$\la = (6,6,4,1,1,1)$ is a partition of length $6$ and weight $19$. Then its conjugate is given by $\la'= (6,3,3,3,2,2)$ and the Ferrers graph of $\la$ is presented by \\

%\input partition19.tex
%Created by jPicEdt 1.x
%Standard LaTeX format (emulated lines)
%Thu Jul 07 15:12:26 CEST 2005
\unitlength 1mm
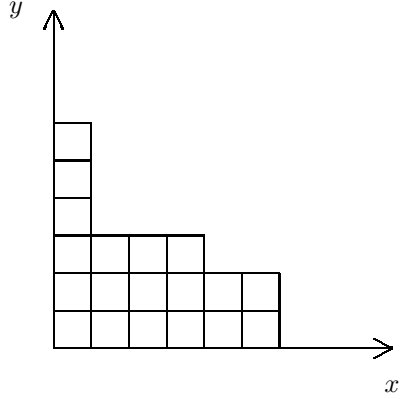
\begin{picture}(85.00,50.00)(0,0)

\linethickness{0.15mm}
%Polygon 0 0(40.00,10.00)(40.00,5.00) 
\put(40.00,5.00){\line(0,1){5.00}}
%End Polygon

\linethickness{0.15mm}
%Polygon 0 0(40.00,5.00)(70.00,5.00) 
\put(40.00,5.00){\line(1,0){30.00}}
%End Polygon

\linethickness{0.15mm}
%Polygon 0 0(70.00,5.00)(70.00,15.00) 
\put(70.00,5.00){\line(0,1){10.00}}
%End Polygon

\linethickness{0.15mm}
%Polygon 0 0(70.00,15.00)(40.00,15.00) 
\put(40.00,15.00){\line(1,0){30.00}}
%End Polygon

\linethickness{0.15mm}
%Polygon 0 0(40.00,15.00)(40.00,10.00) 
\put(40.00,10.00){\line(0,1){5.00}}
%End Polygon

\linethickness{0.15mm}
%Polygon 0 0(40.00,10.00)(70.00,10.00) 
\put(40.00,10.00){\line(1,0){30.00}}
%End Polygon

\linethickness{0.15mm}
%Polygon 0 0(45.00,15.00)(45.00,5.00) 
\put(45.00,5.00){\line(0,1){10.00}}
%End Polygon

\linethickness{0.15mm}
%Polygon 0 0(50.00,15.00)(50.00,5.00) 
\put(50.00,5.00){\line(0,1){10.00}}
%End Polygon

\linethickness{0.15mm}
%Polygon 0 0(55.00,15.00)(55.00,5.00) 
\put(55.00,5.00){\line(0,1){10.00}}
%End Polygon

\linethickness{0.15mm}
%Polygon 0 0(60.00,15.00)(60.00,5.00) 
\put(60.00,5.00){\line(0,1){10.00}}
%End Polygon

\linethickness{0.15mm}
%Polygon 0 0(65.00,15.00)(65.00,5.00) 
\put(65.00,5.00){\line(0,1){10.00}}
%End Polygon

\linethickness{0.15mm}
%Polygon 0 0(40.00,35.00)(40.00,15.00) 
\put(40.00,15.00){\line(0,1){20.00}}
%End Polygon

\linethickness{0.15mm}
%Polygon 0 0(40.00,35.00)(45.00,35.00) 
\put(40.00,35.00){\line(1,0){5.00}}
%End Polygon

\linethickness{0.15mm}
%Polygon 0 0(45.00,35.00)(45.00,15.00) 
\put(45.00,15.00){\line(0,1){20.00}}
%End Polygon

\linethickness{0.15mm}
%Polygon 0 0(40.00,20.00)(60.00,20.00) 
\put(40.00,20.00){\line(1,0){20.00}}
%End Polygon

\linethickness{0.15mm}
%Polygon 0 0(60.00,20.00)(60.00,15.00) 
\put(60.00,15.00){\line(0,1){5.00}}
%End Polygon

\linethickness{0.15mm}
%Polygon 0 0(50.00,20.00)(50.00,15.00) 
\put(50.00,15.00){\line(0,1){5.00}}
%End Polygon

\linethickness{0.15mm}
%Polygon 0 0(55.00,20.00)(55.00,15.00) 
\put(55.00,15.00){\line(0,1){5.00}}
%End Polygon

\linethickness{0.15mm}
%Polygon 0 0(40.00,25.00)(45.00,25.00) 
\put(40.00,25.00){\line(1,0){5.00}}
%End Polygon

\linethickness{0.15mm}
%Polygon 0 0(40.00,30.00)(45.00,30.00) 
\put(40.00,30.00){\line(1,0){5.00}}
%End Polygon

\linethickness{0.15mm}
%Polygon 0 0(40.00,50.00)(40.00,35.00) 
\put(40.00,35.00){\line(0,1){15.00}}
%End Polygon

\linethickness{0.15mm}
%Polygon 0 0(70.00,5.00)(85.00,5.00) 
\put(70.00,5.00){\line(1,0){15.00}}
%End Polygon

\linethickness{0.15mm}
%Polygon 0 0(40.00,50.00)(38.75,47.50) 
\multiput(38.75,47.50)(0.13,0.25){10}{\line(0,1){0.25}}
%End Polygon

\linethickness{0.15mm}
%Polygon 0 0(40.00,50.00)(41.25,47.50) 
\multiput(40.00,50.00)(0.13,-0.25){10}{\line(0,-1){0.25}}
%End Polygon

\linethickness{0.15mm}
%Polygon 0 0(82.50,6.25)(85.00,5.00) 
\multiput(82.50,6.25)(0.25,-0.13){10}{\line(1,0){0.25}}
%End Polygon

\linethickness{0.15mm}
%Polygon 0 0(85.00,5.00)(82.50,3.75) 
\multiput(82.50,3.75)(0.25,0.13){10}{\line(1,0){0.25}}
%End Polygon

\put(35.00,49.38){\makebox(0,0)[cc]{}}

\put(35.00,50.00){\makebox(0,0)[cc]{$y$}}

\put(85.00,0.00){\makebox(0,0)[cc]{$x$}}

\end{picture}
\end{example} 
In the sequel we will omit the axes in Ferrers graphs. For any partition $\la \in \Pscr$ we colour the unit squares of the Ferrers graph $F(\la)$ of $\la$ as follows: a unit square in row $r$ and column $c$ has colour black if $r+c$ is even, and colour white if $r+c$ is odd. The resulting coloured graph is called the {\em chess Ferrers graph} of $\la$. We let $b(\la)$ be the sum of all black unit squares, and $w(\la)$ the sum of all white unit squares. Obviously $b(\la) + w(\la) = n$. 
More formally,
\begin{align*}
b(\la) & = \lceil \frac{\la_{1}}{2} \rceil + \lfloor \frac{\la_{2}}{2} \rfloor + \lceil \frac{\la_{3}}{2} \rceil + \lfloor \frac{\la_{4}}{2} \rfloor + \dots \\
& = \sum_{j}\lceil \frac{\la_{2j+1}}{2} \rceil + \sum_{j} \lfloor \frac{\la_{2j}}{2} \rfloor 
\end{align*}
and
\begin{align*}
w(\la) & = \lfloor \frac{\la_{1}}{2} \rfloor + \lceil \frac{\la_{2}}{2} \rceil + \lfloor \frac{\la_{3}}{2} \rfloor + \lceil \frac{\la_{4}}{2} \rceil + \dots \\
& = \sum_{j}\lfloor \frac{\la_{2j+1}}{2} \rfloor + \sum_{j} \lceil \frac{\la_{2j}}{2} \rceil
\end{align*}
where $\lceil x \rceil$ is the notation for the least integer greater or equal than $x \in \RR$.
\begin{example} \label{example1}
Consider the partition $\la = (6,6,4,1,1,1)$. Then $b(\la) = 9$ and $w(\la) = 10$. The chess Ferrers diagram $F_{\la}$ of $\la$ is given by \\

%\input partition19colour.tex
%Created by jPicEdt 1.x
%Standard LaTeX format (emulated lines)
%Sun Jul 03 18:33:18 CEST 2005
\unitlength 1mm
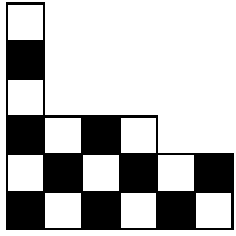
\begin{picture}(70.00,30.00)(0,0)

\linethickness{0.15mm}
%Rectangle(40.00,0.00)(45.00,5.00)  blacken
\put(40.00,0.00){\rule{5.00\unitlength}{5.00\unitlength}}
%End Rectangle

\linethickness{0.15mm}
%Rectangle(45.00,5.00)(50.00,10.00)  blacken
\put(45.00,5.00){\rule{5.00\unitlength}{5.00\unitlength}}
%End Rectangle

\linethickness{0.15mm}
%Rectangle(50.00,0.00)(55.00,5.00)  blacken
\put(50.00,0.00){\rule{5.00\unitlength}{5.00\unitlength}}
%End Rectangle

\linethickness{0.15mm}
%Rectangle(60.00,0.00)(65.00,5.00)  blacken
\put(60.00,0.00){\rule{5.00\unitlength}{5.00\unitlength}}
%End Rectangle

\linethickness{0.15mm}
%Rectangle(55.00,5.00)(60.00,10.00)  blacken
\put(55.00,5.00){\rule{5.00\unitlength}{5.00\unitlength}}
%End Rectangle

\linethickness{0.15mm}
%Rectangle(65.00,5.00)(70.00,10.00)  blacken
\put(65.00,5.00){\rule{5.00\unitlength}{5.00\unitlength}}
%End Rectangle

\linethickness{0.15mm}
%Rectangle(40.00,10.00)(45.00,15.00)  blacken
\put(40.00,10.00){\rule{5.00\unitlength}{5.00\unitlength}}
%End Rectangle

\linethickness{0.15mm}
%Rectangle(50.00,10.00)(55.00,15.00)  blacken
\put(50.00,10.00){\rule{5.00\unitlength}{5.00\unitlength}}
%End Rectangle

\linethickness{0.15mm}
%Rectangle(40.00,20.00)(45.00,25.00)  blacken
\put(40.00,20.00){\rule{5.00\unitlength}{5.00\unitlength}}
%End Rectangle

\linethickness{0.15mm}
%Rectangle(40.00,25.00)(45.00,30.00)  
\put(40.00,25.00){\line(1,0){5.00}}
\put(40.00,25.00){\line(0,1){5.00}}
\put(45.00,25.00){\line(0,1){5.00}}
\put(40.00,30.00){\line(1,0){5.00}}
%End Rectangle

\linethickness{0.15mm}
%Rectangle(40.00,15.00)(45.00,20.00)  
\put(40.00,15.00){\line(1,0){5.00}}
\put(40.00,15.00){\line(0,1){5.00}}
\put(45.00,15.00){\line(0,1){5.00}}
\put(40.00,20.00){\line(1,0){5.00}}
%End Rectangle

\linethickness{0.15mm}
%Rectangle(45.00,15.00)(55.00,15.00)  
\put(45.00,15.00){\line(1,0){10.00}}

\put(45.00,15.00){\line(1,0){10.00}}
%End Rectangle

\linethickness{0.15mm}
%Rectangle(55.00,10.00)(60.00,15.00)  
\put(55.00,10.00){\line(1,0){5.00}}
\put(55.00,10.00){\line(0,1){5.00}}
\put(60.00,10.00){\line(0,1){5.00}}
\put(55.00,15.00){\line(1,0){5.00}}
%End Rectangle

\linethickness{0.15mm}
%Rectangle(40.00,5.00)(45.00,10.00)  
\put(40.00,5.00){\line(1,0){5.00}}
\put(40.00,5.00){\line(0,1){5.00}}
\put(45.00,5.00){\line(0,1){5.00}}
\put(40.00,10.00){\line(1,0){5.00}}
%End Rectangle

\linethickness{0.15mm}
%Rectangle(45.00,0.00)(50.00,5.00)  
\put(45.00,0.00){\line(1,0){5.00}}
\put(45.00,0.00){\line(0,1){5.00}}
\put(50.00,0.00){\line(0,1){5.00}}
\put(45.00,5.00){\line(1,0){5.00}}
%End Rectangle

\linethickness{0.15mm}
%Rectangle(55.00,0.00)(60.00,5.00)  
\put(55.00,0.00){\line(1,0){5.00}}
\put(55.00,0.00){\line(0,1){5.00}}
\put(60.00,0.00){\line(0,1){5.00}}
\put(55.00,5.00){\line(1,0){5.00}}
%End Rectangle

\linethickness{0.15mm}
%Rectangle(60.00,5.00)(65.00,10.00)  
\put(60.00,5.00){\line(1,0){5.00}}
\put(60.00,5.00){\line(0,1){5.00}}
\put(65.00,5.00){\line(0,1){5.00}}
\put(60.00,10.00){\line(1,0){5.00}}
%End Rectangle

\linethickness{0.15mm}
%Rectangle(65.00,0.00)(70.00,5.00)  
\put(65.00,0.00){\line(1,0){5.00}}
\put(65.00,0.00){\line(0,1){5.00}}
\put(70.00,0.00){\line(0,1){5.00}}
\put(65.00,5.00){\line(1,0){5.00}}
%End Rectangle

\linethickness{0.15mm}
%Rectangle(40.00,0.00)(70.00,5.00)  
\put(40.00,0.00){\line(1,0){30.00}}
\put(40.00,0.00){\line(0,1){5.00}}
\put(70.00,0.00){\line(0,1){5.00}}
\put(40.00,5.00){\line(1,0){30.00}}
%End Rectangle

\linethickness{0.15mm}
%Rectangle(40.00,5.00)(70.00,10.00)  
\put(40.00,5.00){\line(1,0){30.00}}
\put(40.00,5.00){\line(0,1){5.00}}
\put(70.00,5.00){\line(0,1){5.00}}
\put(40.00,10.00){\line(1,0){30.00}}
%End Rectangle

\linethickness{0.15mm}
%Rectangle(40.00,10.00)(60.00,15.00)  
\put(40.00,10.00){\line(1,0){20.00}}
\put(40.00,10.00){\line(0,1){5.00}}
\put(60.00,10.00){\line(0,1){5.00}}
\put(40.00,15.00){\line(1,0){20.00}}
%End Rectangle

\linethickness{0.15mm}
%Rectangle(40.00,15.00)(45.00,30.00)  
\put(40.00,15.00){\line(1,0){5.00}}
\put(40.00,15.00){\line(0,1){15.00}}
\put(45.00,15.00){\line(0,1){15.00}}
\put(40.00,30.00){\line(1,0){5.00}}
%End Rectangle

\linethickness{0.15mm}
%Rectangle(45.00,10.00)(45.00,15.00)  

\put(45.00,10.00){\line(0,1){5.00}}
\put(45.00,10.00){\line(0,1){5.00}}

%End Rectangle

\linethickness{0.15mm}
%Rectangle(50.00,10.00)(55.00,15.00)  
\put(50.00,10.00){\line(1,0){5.00}}
\put(50.00,10.00){\line(0,1){5.00}}
\put(55.00,10.00){\line(0,1){5.00}}
\put(50.00,15.00){\line(1,0){5.00}}
%End Rectangle

\linethickness{0.15mm}
%Rectangle(45.00,5.00)(50.00,10.00)  
\put(45.00,5.00){\line(1,0){5.00}}
\put(45.00,5.00){\line(0,1){5.00}}
\put(50.00,5.00){\line(0,1){5.00}}
\put(45.00,10.00){\line(1,0){5.00}}
%End Rectangle

\linethickness{0.15mm}
%Rectangle(55.00,5.00)(60.00,10.00)  
\put(55.00,5.00){\line(1,0){5.00}}
\put(55.00,5.00){\line(0,1){5.00}}
\put(60.00,5.00){\line(0,1){5.00}}
\put(55.00,10.00){\line(1,0){5.00}}
%End Rectangle

\linethickness{0.15mm}
%Rectangle(65.00,5.00)(70.00,10.00)  
\put(65.00,5.00){\line(1,0){5.00}}
\put(65.00,5.00){\line(0,1){5.00}}
\put(70.00,5.00){\line(0,1){5.00}}
\put(65.00,10.00){\line(1,0){5.00}}
%End Rectangle

\linethickness{0.15mm}
%Rectangle(60.00,0.00)(65.00,5.00)  
\put(60.00,0.00){\line(1,0){5.00}}
\put(60.00,0.00){\line(0,1){5.00}}
\put(65.00,0.00){\line(0,1){5.00}}
\put(60.00,5.00){\line(1,0){5.00}}
%End Rectangle

\linethickness{0.15mm}
%Rectangle(50.00,0.00)(55.00,5.00)  
\put(50.00,0.00){\line(1,0){5.00}}
\put(50.00,0.00){\line(0,1){5.00}}
\put(55.00,0.00){\line(0,1){5.00}}
\put(50.00,5.00){\line(1,0){5.00}}
%End Rectangle

\linethickness{0.15mm}
%Rectangle(40.00,0.00)(45.00,5.00)  
\put(40.00,0.00){\line(1,0){5.00}}
\put(40.00,0.00){\line(0,1){5.00}}
\put(45.00,0.00){\line(0,1){5.00}}
\put(40.00,5.00){\line(1,0){5.00}}
%End Rectangle

\linethickness{0.15mm}
%Rectangle(40.00,20.00)(45.00,25.00)  
\put(40.00,20.00){\line(1,0){5.00}}
\put(40.00,20.00){\line(0,1){5.00}}
\put(45.00,20.00){\line(0,1){5.00}}
\put(40.00,25.00){\line(1,0){5.00}}
%End Rectangle

\linethickness{0.15mm}
%Rectangle(40.00,10.00)(45.00,15.00)  
\put(40.00,10.00){\line(1,0){5.00}}
\put(40.00,10.00){\line(0,1){5.00}}
\put(45.00,10.00){\line(0,1){5.00}}
\put(40.00,15.00){\line(1,0){5.00}}
%End Rectangle

\end{picture}
\end{example}

\subsection{From partitions to Castelnuovo functions} \label{2.2}

In the sequel we identify a function $f:\ZZ \r \CC$ with its generating function $f(t)=\sum_n f(n) t^n$. We refer to $f(t)$ as a polynomial or a series depending on whether the support of $f$ is finite or not. \\

A \emph{Castelnuovo function} \cite{Davis} is a finite supported function $s: \NN \r \NN$ such that  
\begin{equation} \label{ref-1.1-1} 
s(0)=1,s(1)=2,\ldots,s(\sigma-1)=\sigma \mbox{ and } s(\sigma-1)\ge 
s(\sigma)\ge s(\sigma+1)\ge \cdots \ge 0
\end{equation} 
for some integer $\sigma \geq 0$. We write $\Sscr$ for the set of all Castelnuovo functions. It is convenient to visualize a Castelnuovo function $s \in \Sscr$ using the graph of the staircase function 
\[ 
F(s): \RR \r \NN: x \mapsto s({\lfloor x \rfloor}) 
\] 
and to divide the area under this graph in unit cases. We will call the result a \emph{Castelnuovo graph} (or {\em Castelnuovo diagram}). The {\em weight} of a Castelnuovo function is the sum of its values, i.e.\ the number of unit squares in the graph. %The {\em height} of $s(t)$ is defined as $\max \{ s_{i} \}$. 
\begin{example} \label{example2}
$s(t) = 1 + 2t + 3t^{2} + 4t^{3} + 5t^{4} + 5t^{5} + 3t^{6} + 2t^{7} + t^{8} + t^{9} + t^{10} + t^{11}$ is a Castelnuovo polynomial of weight $28$. The corresponding Castelnuovo graph is \\

%\input fig29.tex 
%Created by jPicEdt 1.x
%Standard LaTeX format (emulated lines)
%Sun Jul 03 18:27:33 CEST 2005
\unitlength 1mm
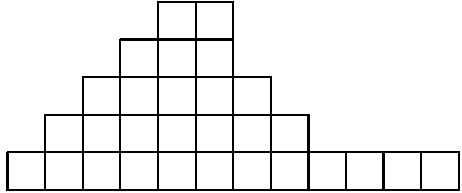
\begin{picture}(90.00,25.00)(0,0)

\linethickness{0.15mm}
%Polygon 0 0(30.00,0.00)(30.00,5.00) 
\put(30.00,0.00){\line(0,1){5.00}}
%End Polygon

\linethickness{0.15mm}
%Polygon 0 0(30.00,5.00)(35.00,5.00) 
\put(30.00,5.00){\line(1,0){5.00}}
%End Polygon

\linethickness{0.15mm}
%Polygon 0 0(35.00,5.00)(35.00,10.00) 
\put(35.00,5.00){\line(0,1){5.00}}
%End Polygon

\linethickness{0.15mm}
%Polygon 0 0(35.00,10.00)(40.00,10.00) 
\put(35.00,10.00){\line(1,0){5.00}}
%End Polygon

\linethickness{0.15mm}
%Polygon 0 0(40.00,10.00)(40.00,15.00) 
\put(40.00,10.00){\line(0,1){5.00}}
%End Polygon

\linethickness{0.15mm}
%Polygon 0 0(40.00,15.00)(45.00,15.00) 
\put(40.00,15.00){\line(1,0){5.00}}
%End Polygon

\linethickness{0.15mm}
%Polygon 0 0(45.00,15.00)(45.00,20.00) 
\put(45.00,15.00){\line(0,1){5.00}}
%End Polygon

\linethickness{0.15mm}
%Polygon 0 0(45.00,20.00)(50.00,20.00) 
\put(45.00,20.00){\line(1,0){5.00}}
%End Polygon

\linethickness{0.15mm}
%Polygon 0 0(50.00,20.00)(50.00,25.00) 
\put(50.00,20.00){\line(0,1){5.00}}
%End Polygon

\linethickness{0.15mm}
%Polygon 0 0(50.00,25.00)(55.00,25.00) 
\put(50.00,25.00){\line(1,0){5.00}}
%End Polygon

\linethickness{0.15mm}
%Polygon 0 0(55.00,25.00)(55.00,0.00) 
\put(55.00,0.00){\line(0,1){25.00}}
%End Polygon

\linethickness{0.15mm}
%Polygon 0 0(55.00,15.00)(60.00,15.00) 
\put(55.00,15.00){\line(1,0){5.00}}
%End Polygon

\linethickness{0.15mm}
%Polygon 0 0(60.00,15.00)(60.00,10.00) 
\put(60.00,10.00){\line(0,1){5.00}}
%End Polygon

\linethickness{0.15mm}
%Polygon 0 0(60.00,10.00)(65.00,10.00) 
\put(60.00,10.00){\line(1,0){5.00}}
%End Polygon

\linethickness{0.15mm}
%Polygon 0 0(65.00,10.00)(65.00,0.00) 
\put(65.00,0.00){\line(0,1){10.00}}
%End Polygon

\linethickness{0.15mm}
%Polygon 0 0(65.00,5.00)(85.00,5.00) 
\put(65.00,5.00){\line(1,0){20.00}}
%End Polygon

\linethickness{0.15mm}
%Polygon 0 0(85.00,5.00)(85.00,0.00) 
\put(85.00,0.00){\line(0,1){5.00}}
%End Polygon

\linethickness{0.15mm}
%Polygon 0 0(85.00,0.00)(30.00,0.00) 
\put(30.00,0.00){\line(1,0){55.00}}
%End Polygon

\linethickness{0.15mm}
%Polygon 0 0(50.00,20.00)(50.00,0.00) 
\put(50.00,0.00){\line(0,1){20.00}}
%End Polygon

\linethickness{0.15mm}
%Polygon 0 0(45.00,15.00)(45.00,0.00) 
\put(45.00,0.00){\line(0,1){15.00}}
%End Polygon

\linethickness{0.15mm}
%Polygon 0 0(40.00,10.00)(40.00,0.00) 
\put(40.00,0.00){\line(0,1){10.00}}
%End Polygon

\linethickness{0.15mm}
%Polygon 0 0(35.00,5.00)(35.00,0.00) 
\put(35.00,0.00){\line(0,1){5.00}}
%End Polygon

\linethickness{0.15mm}
%Polygon 0 0(35.00,5.00)(65.00,5.00) 
\put(35.00,5.00){\line(1,0){30.00}}
%End Polygon

\linethickness{0.15mm}
%Polygon 0 0(40.00,10.00)(60.00,10.00) 
\put(40.00,10.00){\line(1,0){20.00}}
%End Polygon

\linethickness{0.15mm}
%Polygon 0 0(45.00,15.00)(55.00,15.00) 
\put(45.00,15.00){\line(1,0){10.00}}
%End Polygon

\linethickness{0.15mm}
%Polygon 0 0(50.00,20.00)(55.00,20.00) 
\put(50.00,20.00){\line(1,0){5.00}}
%End Polygon

\linethickness{0.15mm}
%Polygon 0 0(60.00,10.00)(60.00,0.00) 
\put(60.00,0.00){\line(0,1){10.00}}
%End Polygon

\linethickness{0.15mm}
%Polygon 0 0(70.00,5.00)(70.00,0.00) 
\put(70.00,0.00){\line(0,1){5.00}}
%End Polygon

\linethickness{0.15mm}
%Polygon 0 0(75.00,5.00)(75.00,0.00) 
\put(75.00,0.00){\line(0,1){5.00}}
%End Polygon

\linethickness{0.15mm}
%Polygon 0 0(80.00,5.00)(80.00,0.00) 
\put(80.00,0.00){\line(0,1){5.00}}
%End Polygon

\linethickness{0.15mm}
%Polygon 0 0(55.00,25.00)(60.00,25.00) 
\put(55.00,25.00){\line(1,0){5.00}}
%End Polygon

\linethickness{0.15mm}
%Polygon 0 0(60.00,25.00)(60.00,15.00) 
\put(60.00,15.00){\line(0,1){10.00}}
%End Polygon

\linethickness{0.15mm}
%Polygon 0 0(55.00,20.00)(60.00,20.00) 
\put(55.00,20.00){\line(1,0){5.00}}
%End Polygon

\linethickness{0.15mm}
%Polygon 0 0(60.00,15.00)(65.00,15.00) 
\put(60.00,15.00){\line(1,0){5.00}}
%End Polygon

\linethickness{0.15mm}
%Polygon 0 0(65.00,15.00)(65.00,10.00) 
\put(65.00,10.00){\line(0,1){5.00}}
%End Polygon

\linethickness{0.15mm}
%Polygon 0 0(65.00,10.00)(70.00,10.00) 
\put(65.00,10.00){\line(1,0){5.00}}
%End Polygon

\linethickness{0.15mm}
%Polygon 0 0(70.00,10.00)(70.00,5.00) 
\put(70.00,5.00){\line(0,1){5.00}}
%End Polygon

\linethickness{0.15mm}
%Polygon 0 0(85.00,5.00)(90.00,5.00) 
\put(85.00,5.00){\line(1,0){5.00}}
%End Polygon

\linethickness{0.15mm}
%Polygon 0 0(90.00,5.00)(90.00,0.00) 
\put(90.00,0.00){\line(0,1){5.00}}
%End Polygon

\linethickness{0.15mm}
%Polygon 0 0(90.00,0.00)(85.00,0.00) 
\put(85.00,0.00){\line(1,0){5.00}}
%End Polygon

\end{picture}
\end{example} 
Given a Castelnuovo function $s$ we colour the unit squares of its Castelnuovo graph $F(s)$ of $s$ as follows: An unit square in column $c$ has colour black if $c$ is even, and colour white if $c$ is odd. Again we agree that the columns are indexed from left to right, and the most left column has index zero. The resulting coloured graph is called the {\em coloured Castelnuovo graph} of $s$. We let $b(s)$ be the sum of all black cases, and $w(s)$ the sum of all white cases. Obviously 
\[
b(s) = \sum_{i}s_{2i}, \quad w(s) = \sum_{i}s_{2i+1}.
\] 
\begin{example} 
For the Castelnuovo polynomial $s(t) = 1 + 2t + 3t^{2} + 4t^{3} + 5t^{4} + 5t^{5} + 3t^{6} + 2t^{7} + t^{8} + t^{9} + t^{10} + t^{11}$ from Example \ref{example2} we have $b(s) = 14$, $w(s) = 15$. The corresponding coloured Castelnuovo graph is given by \\

%\input fig29coloured.tex 
%Created by jPicEdt 1.x
%Standard LaTeX format (emulated lines)
%Thu Jul 07 10:16:28 CEST 2005
\unitlength 1mm
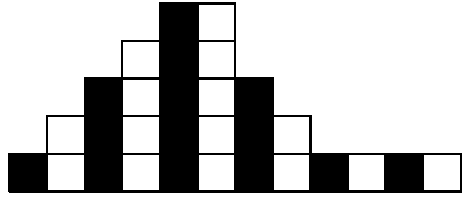
\begin{picture}(90.00,25.00)(0,0)

\linethickness{0.15mm}
%Polygon 0 0(30.00,0.00)(30.00,5.00) 
\put(30.00,0.00){\line(0,1){5.00}}
%End Polygon

\linethickness{0.15mm}
%Polygon 0 0(30.00,5.00)(35.00,5.00) 
\put(30.00,5.00){\line(1,0){5.00}}
%End Polygon

\linethickness{0.15mm}
%Polygon 0 0(35.00,5.00)(35.00,10.00) 
\put(35.00,5.00){\line(0,1){5.00}}
%End Polygon

\linethickness{0.15mm}
%Polygon 0 0(35.00,10.00)(40.00,10.00) 
\put(35.00,10.00){\line(1,0){5.00}}
%End Polygon

\linethickness{0.15mm}
%Polygon 0 0(40.00,10.00)(40.00,15.00) 
\put(40.00,10.00){\line(0,1){5.00}}
%End Polygon

\linethickness{0.15mm}
%Polygon 0 0(40.00,15.00)(45.00,15.00) 
\put(40.00,15.00){\line(1,0){5.00}}
%End Polygon

\linethickness{0.15mm}
%Polygon 0 0(45.00,15.00)(45.00,20.00) 
\put(45.00,15.00){\line(0,1){5.00}}
%End Polygon

\linethickness{0.15mm}
%Polygon 0 0(45.00,20.00)(50.00,20.00) 
\put(45.00,20.00){\line(1,0){5.00}}
%End Polygon

\linethickness{0.15mm}
%Polygon 0 0(50.00,20.00)(50.00,25.00) 
\put(50.00,20.00){\line(0,1){5.00}}
%End Polygon

\linethickness{0.15mm}
%Polygon 0 0(50.00,25.00)(55.00,25.00) 
\put(50.00,25.00){\line(1,0){5.00}}
%End Polygon

\linethickness{0.15mm}
%Polygon 0 0(55.00,25.00)(55.00,0.00) 
\put(55.00,0.00){\line(0,1){25.00}}
%End Polygon

\linethickness{0.15mm}
%Polygon 0 0(55.00,15.00)(60.00,15.00) 
\put(55.00,15.00){\line(1,0){5.00}}
%End Polygon

\linethickness{0.15mm}
%Polygon 0 0(60.00,15.00)(60.00,10.00) 
\put(60.00,10.00){\line(0,1){5.00}}
%End Polygon

\linethickness{0.15mm}
%Polygon 0 0(60.00,10.00)(65.00,10.00) 
\put(60.00,10.00){\line(1,0){5.00}}
%End Polygon

\linethickness{0.15mm}
%Polygon 0 0(65.00,10.00)(65.00,0.00) 
\put(65.00,0.00){\line(0,1){10.00}}
%End Polygon

\linethickness{0.15mm}
%Polygon 0 0(65.00,5.00)(85.00,5.00) 
\put(65.00,5.00){\line(1,0){20.00}}
%End Polygon

\linethickness{0.15mm}
%Polygon 0 0(85.00,5.00)(85.00,0.00) 
\put(85.00,0.00){\line(0,1){5.00}}
%End Polygon

\linethickness{0.15mm}
%Polygon 0 0(85.00,0.00)(30.00,0.00) 
\put(30.00,0.00){\line(1,0){55.00}}
%End Polygon

\linethickness{0.15mm}
%Polygon 0 0(50.00,20.00)(50.00,0.00) 
\put(50.00,0.00){\line(0,1){20.00}}
%End Polygon

\linethickness{0.15mm}
%Polygon 0 0(45.00,15.00)(45.00,0.00) 
\put(45.00,0.00){\line(0,1){15.00}}
%End Polygon

\linethickness{0.15mm}
%Polygon 0 0(40.00,10.00)(40.00,0.00) 
\put(40.00,0.00){\line(0,1){10.00}}
%End Polygon

\linethickness{0.15mm}
%Polygon 0 0(35.00,5.00)(35.00,0.00) 
\put(35.00,0.00){\line(0,1){5.00}}
%End Polygon

\linethickness{0.15mm}
%Polygon 0 0(35.00,5.00)(65.00,5.00) 
\put(35.00,5.00){\line(1,0){30.00}}
%End Polygon

\linethickness{0.15mm}
%Polygon 0 0(40.00,10.00)(60.00,10.00) 
\put(40.00,10.00){\line(1,0){20.00}}
%End Polygon

\linethickness{0.15mm}
%Polygon 0 0(45.00,15.00)(55.00,15.00) 
\put(45.00,15.00){\line(1,0){10.00}}
%End Polygon

\linethickness{0.15mm}
%Polygon 0 0(50.00,20.00)(55.00,20.00) 
\put(50.00,20.00){\line(1,0){5.00}}
%End Polygon

\linethickness{0.15mm}
%Polygon 0 0(60.00,10.00)(60.00,0.00) 
\put(60.00,0.00){\line(0,1){10.00}}
%End Polygon

\linethickness{0.15mm}
%Polygon 0 0(70.00,5.00)(70.00,0.00) 
\put(70.00,0.00){\line(0,1){5.00}}
%End Polygon

\linethickness{0.15mm}
%Polygon 0 0(75.00,5.00)(75.00,0.00) 
\put(75.00,0.00){\line(0,1){5.00}}
%End Polygon

\linethickness{0.15mm}
%Polygon 0 0(80.00,5.00)(80.00,0.00) 
\put(80.00,0.00){\line(0,1){5.00}}
%End Polygon

\linethickness{0.15mm}
%Polygon 0 0(55.00,25.00)(60.00,25.00) 
\put(55.00,25.00){\line(1,0){5.00}}
%End Polygon

\linethickness{0.15mm}
%Polygon 0 0(60.00,25.00)(60.00,15.00) 
\put(60.00,15.00){\line(0,1){10.00}}
%End Polygon

\linethickness{0.15mm}
%Polygon 0 0(55.00,20.00)(60.00,20.00) 
\put(55.00,20.00){\line(1,0){5.00}}
%End Polygon

\linethickness{0.15mm}
%Polygon 0 0(60.00,15.00)(65.00,15.00) 
\put(60.00,15.00){\line(1,0){5.00}}
%End Polygon

\linethickness{0.15mm}
%Polygon 0 0(65.00,15.00)(65.00,10.00) 
\put(65.00,10.00){\line(0,1){5.00}}
%End Polygon

\linethickness{0.15mm}
%Polygon 0 0(65.00,10.00)(70.00,10.00) 
\put(65.00,10.00){\line(1,0){5.00}}
%End Polygon

\linethickness{0.15mm}
%Polygon 0 0(70.00,10.00)(70.00,5.00) 
\put(70.00,5.00){\line(0,1){5.00}}
%End Polygon

\linethickness{0.15mm}
%Polygon 0 0(85.00,5.00)(90.00,5.00) 
\put(85.00,5.00){\line(1,0){5.00}}
%End Polygon

\linethickness{0.15mm}
%Polygon 0 0(90.00,5.00)(90.00,0.00) 
\put(90.00,0.00){\line(0,1){5.00}}
%End Polygon

\linethickness{0.15mm}
%Polygon 0 0(90.00,0.00)(85.00,0.00) 
\put(85.00,0.00){\line(1,0){5.00}}
%End Polygon

\linethickness{0.15mm}
%Rectangle(30.00,0.00)(35.00,5.00)  blacken
\put(30.00,0.00){\rule{5.00\unitlength}{5.00\unitlength}}
%End Rectangle

\linethickness{0.15mm}
%Rectangle(40.00,0.00)(45.00,15.00)  blacken
\put(40.00,0.00){\rule{5.00\unitlength}{15.00\unitlength}}
%End Rectangle

\linethickness{0.15mm}
%Rectangle(50.00,0.00)(55.00,25.00)  blacken
\put(50.00,0.00){\rule{5.00\unitlength}{25.00\unitlength}}
%End Rectangle

\linethickness{0.15mm}
%Rectangle(60.00,0.00)(65.00,15.00)  blacken
\put(60.00,0.00){\rule{5.00\unitlength}{15.00\unitlength}}
%End Rectangle

\linethickness{0.15mm}
%Rectangle(70.00,5.00)(75.00,5.00)  blacken
\put(70.00,5.00){\rule{5.00\unitlength}{0.00\unitlength}}
%End Rectangle

\linethickness{0.15mm}
%Rectangle(70.00,0.00)(75.00,5.00)  blacken
\put(70.00,0.00){\rule{5.00\unitlength}{5.00\unitlength}}
%End Rectangle

\linethickness{0.15mm}
%Rectangle(80.00,0.00)(85.00,5.00)  blacken
\put(80.00,0.00){\rule{5.00\unitlength}{5.00\unitlength}}
%End Rectangle

\end{picture}
\end{example} 
We next describe the relationship between partitions and Castelnuovo functions. For a partition $\la = (\la_{0}, \la_{1}, \dots, \la_{l-l})$ we let $s_{\la} : \NN \r \NN$ be the function defined by 
\[
s_{\la}(m) = \kard \{ j \in \NN \mid j \leq m+1 \text{ and } m+2-j \leq \la_{j} \} 
\]
It is easy to see that $s_{\la}(m)$ is exactly the sum of unit squares which meet the line $D_{m}: y = -x + m$ in the Ferrers graph of $\la$. This corresponds to the interpretation in the introduction.
\begin{example}
Consider the partition $\la = (6,6,4,1,1,1)$ from Example \ref{example1}. We compute
\[
s_{\la}(t) = 1 + 2t + 3t^{2} + 4t^{3} + 4t^{4} + 4t^{5} + t^{6}.
\]
The interpretation for the associated Ferrers graph $F(\la)$ is illustrated for $s_{\la}(1)$ and $s_{\la}(5)$. The line $D_{1}$ meets two unit squares hence $s_{\la}(1) = 2$. Similary $s_{\la}(5) = 4$. \\

%\input compute.tex 
%Created by jPicEdt 1.x
%Standard LaTeX format (emulated lines)
%Thu Jul 07 10:08:28 CEST 2005
\unitlength 1mm
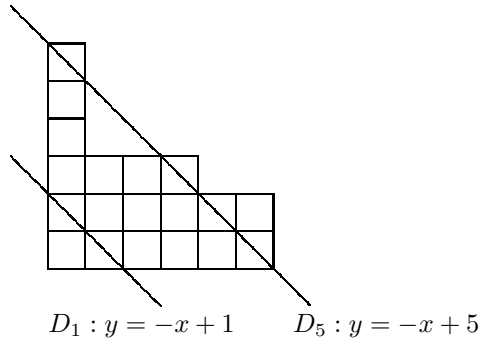
\begin{picture}(85.00,45.00)(0,0)

\linethickness{0.15mm}
%Polygon 0 0(40.00,15.00)(40.00,10.00) 
\put(40.00,10.00){\line(0,1){5.00}}
%End Polygon

\linethickness{0.15mm}
%Polygon 0 0(40.00,10.00)(70.00,10.00) 
\put(40.00,10.00){\line(1,0){30.00}}
%End Polygon

\linethickness{0.15mm}
%Polygon 0 0(70.00,10.00)(70.00,20.00) 
\put(70.00,10.00){\line(0,1){10.00}}
%End Polygon

\linethickness{0.15mm}
%Polygon 0 0(70.00,20.00)(40.00,20.00) 
\put(40.00,20.00){\line(1,0){30.00}}
%End Polygon

\linethickness{0.15mm}
%Polygon 0 0(40.00,20.00)(40.00,15.00) 
\put(40.00,15.00){\line(0,1){5.00}}
%End Polygon

\linethickness{0.15mm}
%Polygon 0 0(40.00,15.00)(70.00,15.00) 
\put(40.00,15.00){\line(1,0){30.00}}
%End Polygon

\linethickness{0.15mm}
%Polygon 0 0(45.00,20.00)(45.00,10.00) 
\put(45.00,10.00){\line(0,1){10.00}}
%End Polygon

\linethickness{0.15mm}
%Polygon 0 0(50.00,20.00)(50.00,10.00) 
\put(50.00,10.00){\line(0,1){10.00}}
%End Polygon

\linethickness{0.15mm}
%Polygon 0 0(55.00,20.00)(55.00,10.00) 
\put(55.00,10.00){\line(0,1){10.00}}
%End Polygon

\linethickness{0.15mm}
%Polygon 0 0(60.00,20.00)(60.00,10.00) 
\put(60.00,10.00){\line(0,1){10.00}}
%End Polygon

\linethickness{0.15mm}
%Polygon 0 0(65.00,20.00)(65.00,10.00) 
\put(65.00,10.00){\line(0,1){10.00}}
%End Polygon

\linethickness{0.15mm}
%Polygon 0 0(40.00,40.00)(40.00,20.00) 
\put(40.00,20.00){\line(0,1){20.00}}
%End Polygon

\linethickness{0.15mm}
%Polygon 0 0(40.00,40.00)(45.00,40.00) 
\put(40.00,40.00){\line(1,0){5.00}}
%End Polygon

\linethickness{0.15mm}
%Polygon 0 0(45.00,40.00)(45.00,20.00) 
\put(45.00,20.00){\line(0,1){20.00}}
%End Polygon

\linethickness{0.15mm}
%Polygon 0 0(40.00,25.00)(60.00,25.00) 
\put(40.00,25.00){\line(1,0){20.00}}
%End Polygon

\linethickness{0.15mm}
%Polygon 0 0(60.00,25.00)(60.00,20.00) 
\put(60.00,20.00){\line(0,1){5.00}}
%End Polygon

\linethickness{0.15mm}
%Polygon 0 0(50.00,25.00)(50.00,20.00) 
\put(50.00,20.00){\line(0,1){5.00}}
%End Polygon

\linethickness{0.15mm}
%Polygon 0 0(55.00,25.00)(55.00,20.00) 
\put(55.00,20.00){\line(0,1){5.00}}
%End Polygon

\linethickness{0.15mm}
%Polygon 0 0(40.00,30.00)(45.00,30.00) 
\put(40.00,30.00){\line(1,0){5.00}}
%End Polygon

\linethickness{0.15mm}
%Polygon 0 0(40.00,35.00)(45.00,35.00) 
\put(40.00,35.00){\line(1,0){5.00}}
%End Polygon

\linethickness{0.15mm}
%Polygon 0 0(35.00,25.00)(55.00,5.00) 
\multiput(35.00,25.00)(0.12,-0.12){167}{\line(1,0){0.12}}
%End Polygon

\put(52.50,2.50){\makebox(0,0)[cc]{$D_{1}: y = -x + 1$}}

\linethickness{0.15mm}
%Polygon 0 0(35.00,45.00)(75.00,5.00) 
\multiput(35.00,45.00)(0.12,-0.12){333}{\line(1,0){0.12}}
%End Polygon

\put(85.00,2.50){\makebox(0,0)[cc]{$D_{5}: y = -x + 5$}}

\end{picture}
\end{example}
The following is immediately clear.
\begin{proposition} \label{correspCastel}
For any partition $\la$ the function $s_{\la}$ is a Castelnuovo function of the same weight. The correspondence $\la \mapsto s_{\la}$ is a surjective map from the set $\Pscr$ of partitions to the set $\Sscr$ of Castelnuovo functions. Furthermore $(b(\la),w(\la)) = (b(s_{\la}),w(s_{\la}))$. 
\end{proposition}
\begin{remark}
As observed in \cite[Remark 1.3]{DV2} follows that the correspondence $\la \mapsto s_{\la}$ restricts to a bijective correspondence between the set $\Dscr$ of partitions in distinct parts and the set $\Sscr$ of Castelnuovo functions. 
\end{remark}

\section{Proof of Theorem A}

\subsection{Proof that the condition in Theorem A is necessary} \label{3.1}

In this subsection we prove that the condition $(b - w)^{2} \leq b$ in Theorem A is necessary. Throughout \S\ref{3.1} $\la \in\Pscr$ is a partition and we denote $(b,w) = (b(\la),w(\la))$.

Consider the map
\begin{align*}
(-)^{\ast} : \ZZ[t] & \r \ZZ[t] \\
f(t) & \mapsto f^{\ast}(t) =
\left\{
\begin{array}{ll}
f(t) - t^{d-1} - t^{d} & \text{ if $f(t) \neq 0$ and $d = \deg f(t) > 0$ } \\
f(t) & \text{ else }
\end{array}
\right. 
%\text{ If $f(t) \neq 0$ and $d = \deg f(t) > 0$, define } f'(t) = f(t) - t^{d-1} - t^{d}
\end{align*}
%for which $f(t) \neq 0$ and $d = \deg f(t) > 0$ we associate the polynomial $f'(t) = f(t) - t^{d-1} - t^{d}$. 
\begin{lemma} \label{maxC}
Assume that $f(t) \neq 0$ is a Castelnuovo polynomial such that $d = \deg f(t) > 0$. If $f^{\ast}(t)$ is not a Castelnuovo polynomial then $f(t)$ is of the form
\[
f(t) = 1 + 2t + 3t^{2} + \dots + (u+1)t^{u}
\]
for some integer $u > 0$.
\end{lemma}
\begin{proof}
Since $f(t)$ is a Castelnuovo polynomial we may write
\[
f(t) = 1 + 2t + 3t^{2} + \dots + (u+1)t^{u} + f_{u+1}t^{u+1} + \dots + f_{v-1}t^{v-1} + f_{v}t^{v}
\]
for some integers $0 \leq u \leq v$ and such that $u+1 \geq f_{u+1} \geq \dots \geq f_{v-1} \geq f_{v} > 0$. It is easy to see that in case $u < v$ then $f^{\ast}(t)$ is a Castelnuovo polynomial. Therefore, if $f^{\ast}(t)$ is not a Castelnuovo polynomial this means $u = v$. This also implies $u > 0$, otherwise $f(t) = 1$ and $\deg f(t) = 0$. Ending the proof.
\end{proof}
Write $s = s_{\la}$ for the Castelnuovo function associated to the partition $\la$. Proposition \ref{correspCastel} implies $(b,w) = (b(s),w(s))$. %If $s(t) = 0$ or $s(t) = 1$ then \eqref{ineq} is checked, thus we may assume that $\deg s(t) > 0$. 
We put 
\[
s_{0}(t) = s(t), s_{1}(t) = s^{\ast}(t), s_{2}(t) = s^{\ast \ast}(t), \dots
\] 
Either $s_{k}$ is a Castelnuovo function for all integers $k \in \NN$, or not. We will treat these two cases seperately.
\begin{case}
$s_{k}$ is a Castelnuovo function for all integers $k \in \NN$.
\end{case}
It is clear that $s_{k} = s_{k+1}$ implies $s_{k+1} = s_{k+2}$ for all integers $k \in \NN$.
Define
\[
l = \max \{ k \in \NN \mid s_{k} \neq s_{k+1} \} + 1
\]
Then $s_{0} \neq s_{1} \neq \dots \neq s_{l-1} \neq s_{l} = s_{l+1} = s_{l+2} = \dots$.
By definition of the map $(-)^{\ast}$ and the fact that $s_{k}$ is a Castelnuovo function we deduce either $s_{l}(t) = 1$ or $s_{l}(t) = 0$. Since for all $k \in \NN$
\[
(b,w) = (b(s),w(s)) = (b(s_{k}) + k,w(s_{k})+k)  
\]
we either have $(b,w) = (l,l)$ or $(b,w) = (l+1,l)$, for which \eqref{ineq} is easily checked.  
\begin{case}
There exists an integer $k$ such that $s_{k}$ is {\em not} a Castelnuovo function.
\end{case}
Put
\[
l = \max \{ k \in \NN \mid s_{k} \text{ is a Castelnuovo function} \}
\] 
This definition makes sense because $s = s_{0}$ is a Castelnuovo function. Lemma \ref{maxC} implies that $s_{l}(t)$ is of the form
\[
s_{l}(t) = 1 + 2t + 3t^{2} + \dots + (u+1)t^{u}
\]
for some integer $u > 0$. One easily computes
\begin{eqnarray} \label{u}
(b(s_{l}),w(s_{l})) = 
\left\{
\begin{array}{ll}
\left((u+2)^{2}/4, u(u+2)/4 \right) & \text{ if $u$ is even } \\
\left((u+1)^{2}/4, (u+1)(u+3)/4 \right) & \text{ if $u$ is odd }
\end{array}
\right.
\end{eqnarray}
and combining with $(b,w) = (b(s),w(s)) = (b(s_{l})+l,w(s_{l})+l)$ we find 
\[
\frac{1}{b} - \left(1-\frac{w}{b} \right)^{2} = \frac{l}{b} \geq 0
\]
which completes the proof.

\subsection{Proof that the condition in Theorem A is sufficient}

Let $b,w \in \NN$ be positive integers such that \eqref{ineq} holds. If $b = 0$ then it follows that $w = 0$, and it is clear that for the empty partition $\la = (\,\,)$ we have $(b,w) = (0,0) = (b(\la),w(\la))$. Hence we may assume $b > 0$. Let
\[
l = \max \{ j \in \NN \mid \sum_{i = 0}^{j}(2i+1) \leq b \text{ and } \sum_{i = 0}^{j}2i \leq w \}
\]
It is clear that there exist positive integers $b',w' \in \NN$ for which either Case 1 or Case 2 is true
\begin{align*}
\text{ Case 1: } & 
\left\{
\begin{array}{l}
b = 1 + 3 + 5 + \dots + (2l-1) + b' \\
w = 2 + 4 + 6 + \dots + 2l + w'
\end{array}
\right.
\text{ where } b' < 2l+1
\end{align*}
\begin{align*}
\text{ Case 2: } & 
\left\{
\begin{array}{l}
b = 1 + 3 + 5 + \dots + (2l+1) + b' \\
w = 2 + 4 + 6 + \dots + 2l + w'
\end{array}
\right.
\text{ where } w' < 2l+2
\end{align*}
\begin{lemma} \label{lemma}
Let $b,w \in \NN$ such that \eqref{ineq} holds, i.e. $(b-w)^{2} \leq b$. Consider the associated integers $l, b', w' \in \NN$ as defined above. We have
\begin{enumerate}
\item
If Case 1 is true then $w' \leq b'$, and  
\item
if Case 2 is true then $b' \leq w'$.   
\end{enumerate}
\end{lemma}
\begin{proof}
\begin{enumerate}
\item
First assume Case 1 is true. Then 
\begin{align*}
\left\{
\begin{array}{l}
b = l^{2} + b' \\
w = l(l+1) + w'
\end{array}
\right.
\end{align*}
From the inequality \eqref{ineq} we find $0 \leq b - (b - w)^{2}$ hence 
\begin{align*}
0 & \leq (l^{2} + b') - \left( l^{2} + b' - l(l+1) - w' \right)^{2} \\
& = l^{2} + b' - (b' - w' - l)^{2} \\
%& = l^{2} + b' - (b' - w')^{2} - l^{2} + 2(b'-w')l \\
& = b' - (b' - w')^{2} + 2(b'-w')l 
\end{align*}
Assume by contradiction $w' > b'$ i.e. $b' - w' \leq -1$. Then we further deduce
\begin{align*}
0 & \leq b' - (b' - w')^{2} + 2(b'-w')l \\
& < b' - (b' - w')^{2} - 2l \\
& \leq -(b' - w')^{2}
\end{align*}
where we have used $b' \leq 2l$. We conclude $0 < -(b' - w')^{2}$, clearly a contradiction. Hence $w' \leq b'$.
\item
Second, assume Case 2 is true. We now have 
\begin{align*}
b & = (l+1)^{2} + b' \\
w & = l(l+1) + w'
\end{align*}
and $0 \leq b - (b - w)^{2}$ leads to
\begin{align*}
0 & \leq \left( (l+1)^{2} + b'\right) - \left( (l+1)^{2} + b' - l(l+1) - w' \right)^{2} \\
& = (l+1)^{2} + b' - \left( (b' - w') + (l + 1) \right)^{2} \\
%& = (l+1)^{2} + b' - (b' - w')^{2} - (l + 1)^{2} - 2(b' - w')(l+1) \\
& = b' - (b' - w')^{2} - 2(b' - w')(l+1) 
\end{align*}
Assume by contradiction $w' < b'$. This means $1 \leq b' - w'$ and also $(b' - w') \leq (b' - w')^{2}$. Invoking these inequalities we further deduce
\begin{align*}
0 & \leq b' - (b' - w')^{2} - 2(b' - w')(l+1) \\
& \leq b' - (b' - w') - 2(b' - w')(l+1) \\
& \leq b' - (b' - w') - 2(l+1)
\end{align*}
and therefore
\[
2l+2 \leq w'
\]
which contradicts $w' < 2l+2$. We conclude $w' \geq b'$, which proves the lemma.
\end{enumerate}
\end{proof}
We now put
\begin{eqnarray*}
s(t) = 
\left\{
\begin{array}{ll}
1 + 2t + 3t^{2} + \dots + (2l-1)t^{2l-2} + (2l)t^{2l-1} + b't^{2l} + w't^{2l+1} & \text{ if Case 1} \\
1 + 2t + 3t^{2} + \dots + (2l)t^{2l-1} + (2l+1)t^{2l} + w't^{2l+1} + b't^{2l+2} & \text{ if Case 2} 
\end{array}
\right.
\end{eqnarray*}
As a consequence of Lemma \ref{lemma} we have that $s(t)$ is a Castelnuovo polynomial for which $(b(\la),w(\la)) = (b,w)$. By Proposition \ref{correspCastel} there exists a partition (in distinct parts) $\la$ for which $(b(\la),w(\la)) = (b,w)$. This proves that the condition \eqref{ineq} in Theorem A is sufficient.

\section{Proof of Theorem B}

In this section we prove Theorem B. First let $\la \in \Pscr$ be any partition. As shown in Section \ref{3.1} there exists integers $k,l$ for which $(b(\la),w(\la))$ is either equal to 
\begin{itemize}
\item
$(l,l)$, or
\item
$(l+1,l)$, or
\item
$\left((k+1)^{2} + l, k(k+1)+ l \right)$ \quad (put $k = u/2$ in \eqref{u} if $u$ is even ), or
\item
$\left(k^{2} + l, k(k+1)+ l \right)$ \quad (put $k = (u+1)/2$ in \eqref{u} if $u$ is odd ).
\end{itemize}
Hence there exist positive integers $k,l \in \NN$ such that either
\[
(b,w) = \left( (k+1)^{2} + l ,k(k+1) + l \right)
\]
or
\[
(b,w) = \left( k^{2} + l, k(k+1) + l \right)
\]
Conversely, let $k,l \in \NN$. Putting  
\[
(b,w) = \left( (k+1)^{2} + l ,k(k+1) + l \right)
\]
it is easy to verify that $b - (b-w)^{2} = l$. Hence \eqref{ineq} holds. By Theorem A there exists a partition $\la$ such that $(b(\la),w(\la)) = (b,w)$. Similar treatment if we put $(b,w) = \left( k^{2} + l, k(k+1) + l \right)$. This ends the proof of Theorem B.

\section{A reformulation}

In this final part we make the connection with Problem 10 of \cite{competition}. For convenience for the reader we recall the quastion as it was stated in \cite{competition}.
\setcounter{problem}{9}
\begin{problem}
Let $n$ be a positive integer. Let $a_{1},a_{2},\dots,a_{m}$ be a partition of $n$. Represent this partition as a left-justified array of boxes, with $a_{1}$ boxes in the first row, $a_{2}$ in the second, and so on, and label the boxes with $1$ and $-1$ in a chess-board pattern, starting with a $1$ in the top-left corner. Let $c$ be the sum of these labels. For instance, if $n = 11$ and the partition is $4,3,3,1$ then $c = -1$, as one sees by summing the labels in the diagram: \\

%\input competition.tex \\
%Created by jPicEdt 1.x
%Standard LaTeX format (emulated lines)
%Wed Jul 13 22:47:57 CEST 2005
\unitlength 1mm
\begin{picture}(65.00,20.00)(0,0)

\linethickness{0.15mm}
%Rectangle(45.00,15.00)(50.00,20.00)  
\put(45.00,15.00){\line(1,0){5.00}}
\put(45.00,15.00){\line(0,1){5.00}}
\put(50.00,15.00){\line(0,1){5.00}}
\put(45.00,20.00){\line(1,0){5.00}}
%End Rectangle

\linethickness{0.15mm}
%Rectangle(50.00,15.00)(55.00,20.00)  
\put(50.00,15.00){\line(1,0){5.00}}
\put(50.00,15.00){\line(0,1){5.00}}
\put(55.00,15.00){\line(0,1){5.00}}
\put(50.00,20.00){\line(1,0){5.00}}
%End Rectangle

\linethickness{0.15mm}
%Rectangle(55.00,15.00)(60.00,20.00)  
\put(55.00,15.00){\line(1,0){5.00}}
\put(55.00,15.00){\line(0,1){5.00}}
\put(60.00,15.00){\line(0,1){5.00}}
\put(55.00,20.00){\line(1,0){5.00}}
%End Rectangle

\linethickness{0.15mm}
%Rectangle(60.00,15.00)(65.00,20.00)  
\put(60.00,15.00){\line(1,0){5.00}}
\put(60.00,15.00){\line(0,1){5.00}}
\put(65.00,15.00){\line(0,1){5.00}}
\put(60.00,20.00){\line(1,0){5.00}}
%End Rectangle

\linethickness{0.15mm}
%Rectangle(45.00,10.00)(50.00,15.00)  
\put(45.00,10.00){\line(1,0){5.00}}
\put(45.00,10.00){\line(0,1){5.00}}
\put(50.00,10.00){\line(0,1){5.00}}
\put(45.00,15.00){\line(1,0){5.00}}
%End Rectangle

\linethickness{0.15mm}
%Rectangle(50.00,10.00)(55.00,15.00)  
\put(50.00,10.00){\line(1,0){5.00}}
\put(50.00,10.00){\line(0,1){5.00}}
\put(55.00,10.00){\line(0,1){5.00}}
\put(50.00,15.00){\line(1,0){5.00}}
%End Rectangle

\linethickness{0.15mm}
%Rectangle(55.00,10.00)(60.00,15.00)  
\put(55.00,10.00){\line(1,0){5.00}}
\put(55.00,10.00){\line(0,1){5.00}}
\put(60.00,10.00){\line(0,1){5.00}}
\put(55.00,15.00){\line(1,0){5.00}}
%End Rectangle

\linethickness{0.15mm}
%Rectangle(45.00,5.00)(50.00,10.00)  
\put(45.00,5.00){\line(1,0){5.00}}
\put(45.00,5.00){\line(0,1){5.00}}
\put(50.00,5.00){\line(0,1){5.00}}
\put(45.00,10.00){\line(1,0){5.00}}
%End Rectangle

\linethickness{0.15mm}
%Rectangle(50.00,5.00)(55.00,10.00)  
\put(50.00,5.00){\line(1,0){5.00}}
\put(50.00,5.00){\line(0,1){5.00}}
\put(55.00,5.00){\line(0,1){5.00}}
\put(50.00,10.00){\line(1,0){5.00}}
%End Rectangle

\linethickness{0.15mm}
%Rectangle(55.00,5.00)(60.00,10.00)  
\put(55.00,5.00){\line(1,0){5.00}}
\put(55.00,5.00){\line(0,1){5.00}}
\put(60.00,5.00){\line(0,1){5.00}}
\put(55.00,10.00){\line(1,0){5.00}}
%End Rectangle

\linethickness{0.15mm}
%Rectangle(45.00,0.00)(50.00,5.00)  
\put(45.00,0.00){\line(1,0){5.00}}
\put(45.00,0.00){\line(0,1){5.00}}
\put(50.00,0.00){\line(0,1){5.00}}
\put(45.00,5.00){\line(1,0){5.00}}
%End Rectangle

\put(47.50,17.50){\makebox(0,0)[cc]{$1$}}

\put(52.50,17.50){\makebox(0,0)[cc]{$-1$}}

\put(57.50,17.50){\makebox(0,0)[cc]{$1$}}

\put(62.50,17.50){\makebox(0,0)[cc]{$-1$}}

\put(47.50,12.50){\makebox(0,0)[cc]{$-1$}}

\put(52.50,12.50){\makebox(0,0)[cc]{$1$}}

\put(57.50,12.50){\makebox(0,0)[cc]{$-1$}}

\put(47.50,7.50){\makebox(0,0)[cc]{$1$}}

\put(52.50,7.50){\makebox(0,0)[cc]{$-1$}}

\put(57.50,7.50){\makebox(0,0)[cc]{$1$}}

\put(47.50,2.50){\makebox(0,0)[cc]{$-1$}}

\end{picture} \\

Prove that $n \geq c(2c-1)$, and determine when equality occurs.
\end{problem}
Let us now indicate how we use Theorem A and Theorem B to solve Problem 10. Write $\la = (a_{1},a_{2},\dots,a_{m})$, and put $(n(\la),c(\la)) = (n,c)$ and $(b,w) = (b(\la),w(\la))$. It is clear that $n = b+w$, $c = b - w$. Hence $b = (n+c)/2$, $w = (n-c)/2$ and it follows that $n + c$ and $n - c$ are even, i.e. $n$ and $c$ have the same parity (either $n$ and $c$ are both even, or they are both odd). Further inequality \eqref{ineq} is equivalent with 
\begin{align*}
(b-w)^{2} \leq b & \Leftrightarrow \left( \frac{n+c}{2} - \frac{n-c}{2} \right) \leq \frac{n+c}{2} \\
& \Leftrightarrow 2c^{2} \leq n+c \\
& \Leftrightarrow c(2c-1) \leq n
\end{align*}
Hence Theorem A implies that $c(2c-1) \leq n$. Conversely, given any $(n,c) \in \NN \times \ZZ$ of the same parity for which $c(2c-1) \leq n$ holds, we see that by putting $b = (n+c)/2$, $w = (n-c)/2$ that \eqref{ineq} holds, hence Theorem A implies that there exists a partition $\la$ such that $(n(\la),c(\la)) = (n,c)$. \\

To see when equality in $c(2c-1) \leq n$ occurs, we may invoke Theorem B: The appearing integers $b,w$ are of the form
\[
(b,w) = \left( (k+1)^{2} + l ,k(k+1) + l \right) \text{ or } (b,w) = \left( k^{2} + l, k(k+1) + l \right)
\]
for some $k,l \in \NN$, and conversely for any $(b,w)$ of this form there exists a partition $\la$ for which $(b,w) = (b(\la),w(\la))$. By replacing $b = (n+c)/2$, $w = (n-c)/2$ we find
\begin{equation} \label{forms}
(n,c) = (2k^{2} + k + 2l,-k) \text{ or } (n,c) = (2k^{2}+3k+1+2l,k+1) 
\end{equation}
for some $k,l \in \NN$, and conversely for any $(n,c)$ of this form there exists a partition $\la$ for which $(n,c) = (n(\la),c(\la))$. Hence for any $c \in \ZZ$ the appearing $n \in \ZZ$ for which \eqref{forms} holds are
\[
n = c(2c - 1) + 2l, \quad l \in \NN.
\]
Note that it follows that $n \in \NN$. Hence equality in $c(2c - 1) \leq n$ occurs if and only if $l = 0$. Using the resuls of section \ref{3.1} we find that $n = c(2c - 1)$ if and only if the associated Castelnuovo function is of the ''maximal" form from the introduction, i.e. the partition is of the form $\la = (m,m-1,\dots,2,1)$ for some $m \in \NN$. We have proved
\setcounter{solution}{9}
\begin{solution}[To Problem 10]
Let $(n,c) \in \NN \times \ZZ$. Then there exists a partition $\la$ such that $(n(\la),c(\la)) = (n,c)$ if and only if 
\[
n,c \text{ have the same parity and } c(2c-1) \leq n 
\]
In this case, $n = c(2c - 1) + 2l$ for some $l \in \NN$. For any partition $\la$ we have $c\left(2c-1\right) = n$ if and only if $\la = (m,m-1,\dots,2,1)$ for some $m \in \NN$. \\
Furthermore the same statement holds if we restrict ourselves to partitions in distinct parts.
\end{solution}
\begin{remark}
The reader will notice that the presented solution of Problem 10 is different from the one presented in \cite[Problem 10]{competition}. Our version is somewhat longer, however the description is more detailed as we alse give the necessary conditions for $(n,c)$ to correspond to a partition. As a consequence, for any partition $\la$ the difference of $n$ and $c(2c - 1)$ is always even. 
\end{remark}

\ifx\undefined\bysame 
\newcommand{\bysame}{\leavevmode\hbox to3em{\hrulefill}\,} 
\fi

\end{document}